\documentclass[article]{article}

\usepackage[utf8]{inputenc}
\usepackage[ngerman,english]{babel}
\usepackage{amsmath}
\usepackage{amsthm}
\usepackage{bm}
\usepackage{amssymb}
\usepackage{graphicx}
\usepackage{multirow}
\usepackage{rotating}
\usepackage{listings}
\usepackage{xcolor} % for setting colors
\usepackage{wrapfig}
\usepackage{lipsum}  % generates filler text
\usepackage{tikz}
\usetikzlibrary{arrows,matrix,positioning}
\usepackage{float}
\usepackage{booktabs}
\usepackage{parskip}
\usepackage{enumitem}
\usepackage{tcolorbox }
\graphicspath{{data/}}
\usepackage[font=small,labelfont=bf]{caption}
\usepackage{lmodern}
\usepackage{subcaption}
\usepackage[left=3cm,right=3cm,bottom=3cm]{geometry}
\usepackage{tikz}
\usepackage{url}
\usetikzlibrary{arrows, decorations.markings}
\usepackage{xcolor}
\usepackage{pdfpages}
\usepackage{siunitx}
\usepackage{color,soul}

\usepackage{todonotes}
\usepackage{svg}
\usepackage{epstopdf}
\usepackage[percent]{overpic}

\usepackage{authblk}

\usepackage{mathtools,algorithm}
\usepackage{algpseudocode}
\algnewcommand\And{\textbf{and} }
\algnewcommand\Or{\textbf{or} }

\algnewcommand{\Inputs}[1]{%
  \State \textbf{Inputs:}
  \Statex \hspace*{\algorithmicindent}\parbox[t]{.8\linewidth}{\raggedright #1}
}
\usepackage{amsthm}
\DeclareMathOperator*{\argmin}{arg\,min}
\theoremstyle{definition}

\newtheorem{theorem}{Theorem}[section]

\newtheorem{corollary}{Corollary}[theorem]

\newtheorem*{remark}{Remark}

\usepackage{siunitx}
\usepackage{mathtools}

\usepackage{algorithm}% http://ctan.org/pkg/algorithms
\usepackage{algpseudocode}% http://ctan.org/pkg/algorithmicx
\algdef{SE}[DOWHILE]{Do}{doWhile}{\algorithmicdo}[1]{\algorithmicwhile\ #1}%
\newcommand{\R}{\mathbb{R}}
\newcommand{\N}{\mathbb{N}}
\newcommand{\card}{\mathop\mathrm{card}}
\newcommand{\dom}{\mathop\mathrm{dom}}

\title{Adaptive Gaussian Process Regression for Efficient Building of Surrogate Models in Inverse Problems}
\author[1]{Phillip Semler}
\author[1]{Martin Weiser}
\affil[1]{Zuse Institute Berlin, Takustr. 7, 14195 Berlin, Germany}
\date{\today}

\begin{document}
\maketitle

\begin{abstract}
In a task where many similar inverse problems must be solved, evaluating costly simulations is impractical. Therefore, replacing the model $y$ with a surrogate model $y_s$ that can be evaluated quickly leads to a significant speedup. The approximation quality of the surrogate model depends strongly on the number, position, and accuracy of the sample points. With an additional finite computational budget, this leads to a problem of (computer) experimental design. In contrast to the selection of sample points, the trade-off between accuracy and effort has hardly been studied systematically. We therefore propose an adaptive algorithm to find an optimal design in terms of position and accuracy. Pursuing a sequential design by incrementally appending the computational budget leads to a convex and constrained optimization problem. As a surrogate, we construct a Gaussian process regression model. We measure the global approximation error in terms of its impact on the accuracy of the identified parameter and aim for a uniform absolute tolerance, assuming that $y_s$ is computed by finite element calculations. A priori error estimates and a coarse estimate of computational effort relate the expected improvement of the surrogate model error to computational effort, resulting in the most efficient combination of sample point and evaluation tolerance. We also allow for improving the accuracy of already existing sample points by continuing previously truncated finite element solution procedures.
\end{abstract}

\section{Introduction}

Any physical measurement or physical model can be formally described by the functional relation $g = y(p)$, where the model $y(p)$ maps parameters and inputs $p$ to observable quantities $g$. Inferring model parameters $p$ from measurements $g$ in order to gain information on the system's state is known as an inverse problem, in contrast to the forward problem of computing or predicting measurement data from parameters, i.e. evaluating $y$. Inverse problems are usually addressed by computing point estimates $p_*$ by solving a minimization problem such as $p_* = \mathrm{arg\,min}_p \|y(p)-g\|$ or by sampling the posterior probability distribution of the parameters conditioned on the available measurement data, see, e.g.,~\cite{Aster, EngelHankeNeubauer1996, KaipioSomersalo2005, Neto, Tarantola}.

In many applications, the model $y$ is not given analytically, but only in form of a complex numerical procedure such as solving a partial differential equation with finite elements. Then, both solving an optimisation problem for computing a point estimate and sampling the posterior distribution require a significant number of model evaluations, rendering the inverse problem a computationally demanding task. In online and real-time applications such as quality control and nondestructive testing, where many similarly structured inverse problems need to be solved in short time, the computational effort of evaluating the model can be prohibitive, such that the original procedure for evaluating $y$ needs to be replaced by a much faster alternative. In an online-offline splitting, first a replacement for $y$ is constructed in an offline phase and used later in the online phase for actually solving the inverse problems. 

Besides model reduction, where finite element discretizations of PDEs arising in $y$ are replaced by a lower dimensional and problem-adapted discretization~\cite{MOR}, surrogate models are used for a direct approximation of $y$ itself, as long as both parameters and measurements are low-dimensional. Popular types of surrogate models or response surfaces include polynomials~\cite{Owen}, sparse grids~\cite{Xiang}, tensor trains~\cite{Griebel,Bigoni}, artificial neural networks~\cite{Holena,Raissi}, and Gaussian process regression (GPR)~\cite{Schneider,Zaytsev}. All of these surrogate models interpolate or approximate $y$ based on function values $y(p_i)$ at some parameter sample points $p_i$, which form the training data. The approximation quality of the resulting surrogate model depends highly on the number and position of these sample points. While in principle an arbitrary amount of training data can be generated simply by numerically evaluating $y$ in the offline phase, building a good surrogate model can be computationally very expensive if the required number of sample points $p_i$ is large. Consequently, design of computer experiment strategies for a near optimal sample point selection have been proposed in order to reduce the number of sample points necessary for achieving a desired approximation accuracy, in particular for the case of analytically well-understood GPR~\cite{Rasmussen}.
A priori static sample point sets defined by factorial designs~\cite{Giunta} or space-filling designs~\cite{Queipo} are complemented by adaptive designs~\cite{Lehmensiek,Sugiyama,Crombecq,Joseph}. Here, the selection is in general based on a readily available pointwise estimate of the surrogate approximation error, and including the parameter point that maximizes a certain acquisition function into the sample set. Acquisition functions such as confidence bounds~\cite{Srinivas}, expected improvement~\cite{Mockus}, probability of improvement~\cite{Kushner}, Thompson sampling~\cite{Thompson}, entropy search~\cite{Hennig}, and knowledge gradient \cite{Wu} are primarily used in Bayesian optimization. For offline created surrogate models, a uniform approximation accuracy is often desired, and thus often a maximizer of the error estimator is chosen as next sample point. 

When computing training data at sample points $p_i$ by a numerical procedure such as a finite element solver, the resulting evaluations of $y$ are always inexact due to discretization and truncation errors. Again, in principle a highly accurate evaluation is possible, but comes with a corresponding computational effort. This accuracy-effort trade-off is a second optimization dimension for the design of computer experiments. In contrast to the selection of parameter sample points, however, it has so far barely been investigated systematically. The use of two different model accuracies, a high-fidelity and a low-fidelity model, has been proposed in~\cite{Nitzler}, and attempts into taking evaluation accuracy into account have been made in~\cite{SagnolHegeWeiser2016}.

Here, we consider constructing a GPR surrogate model offline for the purpose of online parameter identification over a bounded domain, and devise a greedy-type strategy for the sequential computation of training data by selecting simultaneously the next sample position and evaluation accuracy. In a goal-oriented approach~\cite{BeckerRannacher2001}, we measure the surrogate approximation error by its impact on the accuracy of the identified parameter, and aim at a uniform absolute tolerance or, if this cannot be achieved, at least a uniform bounded deterioration with respect to the exact model.
We focus on the case where $y$ is dominated by an error-controlled finite element computation based on adaptive grid refinement. Standard a priori error estimates and a coarse estimate of the computational work incurred by a prescribed tolerance then allows relating the expected improvement of the surrogate model approximation error to the computational work spent, and thus selecting the best, i.e. most efficient, combination of sample point and evaluation tolerance for computing the next training data. In that context, we also allow for improving the accuracy of already existing sample points by continuing previously truncated finite element solution procedures.

The remainder of the paper is organized as follows:
Section~\ref{sec:identification} describes the inversion method used. We start by defining and describing the forward problem for generating the training data. Then, the inverse problem for parameter reconstruction is presented in the notational frame of Bayesian inversion.
Section~\ref{sec:active-learning} presents the GPR surrogate model, describes the accuracy and work models, and works out the greedy heuristic for designing the training data simulations.
Effectiveness and efficiency of the adaptive surrogate model training are investigated at numerical examples in Section~\ref{sec:experiments}.

\section{Surrogate-based parameter identification} \label{sec:identification}

In the following, we set the frame for parameter identification by Tikhonov regularization, equivalent to computing a maximum posterior point estimate in a Bayesian perspective. The reconstruction  depends on a forward problem describing the dependence of the model output $y(p)$ on the parameters $p$, an optimization procedure for the maximum posterior problem, and provides uncertainty quantification via Laplace's method.

\subsection{Forward problem}
The forward problem maps the parameter vector $p\in\R^d$ to the model output $y(p)\in\R^m$ representing the physically measurable quantities. Possible pairings are, for example, structural parameters and occurring strains from structural mechanical problems, geometrical parameters of micro structures and light intensities from scatterometry, or thermal parameters and temperature values from heat conduction, to name a few. We assume the model $y$ to be twice continuously differentiable.

We assume the relevant parameter space $\mathcal{X}\subset \R^d$ to be closed and bounded, usually given by simple bounds $\underbar{p} \le p \le \bar{p}$ or in terms of inequality constraints $g_i(p)\le 0$, $i=1,\dots,r$. On one hand, this allows excluding non-physical parameters explicitly from the reconstruction, i.e. negative densities or absolute temperatures, and, on the other hand, renders a faithful approximation of the model $y$ by some surrogate model on $\mathcal{X}$ practically feasible.

We assume that evaluating $y(p)$ exactly is computationally expensive or even infeasible, but that approximations $y_\epsilon(p)$ with an error $\|y_\epsilon(p)-y(p)\|_{\mathcal{Y}} \le \epsilon$ in some problem-dependent norm $\|\cdot\|_{\mathcal{Y}} $ can be obtained for any prescribed tolerance $\epsilon>0$ with finite computational cost, usually growing for $\epsilon\to 0$. This is the situation we face, i.e., if the model output is defined in terms of solutions of ordinary or partial differential equations, where  the finite element discretization error depends on the mesh width and therefore problem size.

\subsection{Inverse problem} \label{sec:inverse-problem}
In a Bayesian framework, the posterior,  i.e. the conditional probability of parameters $p$ given the measurements $y^m$, is  the product
\begin{equation}\label{eq:postdist}
\pi(p | y^m) \propto  \pi(y^m | p)\pi(p)
\end{equation}
of likelihood $\pi(y^m|p)$ connected to measurement errors and prior probability $\pi(p)$ expressing the a priori knowledge on parameters.  For simplicity of notation we restrict the attention to Gaussian likelihood
\begin{equation}\label{eq:likelihood}
\pi(y^m | p) \sim \exp{\left ( -\frac{1}{2} ( y^m-y(p)  )^T\Sigma^{-1}_l  (y^m-y(p))  \right )} \quad
\end{equation}
and prior
\begin{equation}\label{eq:prior}
\pi(p) \sim \exp\left ( -\frac{1}{2}( p-p^0)^T\Sigma^{-1}_p  ( p-p^0) \right ) 
\end{equation}
with covariance $\Sigma_l\in\R^{m\times m}$ and $\Sigma_p\in\R^{d\times d}$, respectively. The model output $y(p)$ provides the likelihood mean, whereas the prior mean $p^0$ is constant.

An expensive evaluation of the model $y(p)$ renders a sampling of the posterior by Markov Chain Monte Carlo methods~\cite{MonteCarlo} unattractive due to the generally large number of likelihood evaluations required. We therefore consider the reconstruction of parameters via a maximum posterior point estimate $p_{\rm MAP} = p(y^m) := \mathop\mathrm{arg\; max}_{p\in\mathcal{X}} \pi(p|y^m)$, which amounts to minimizing
\begin{equation}\label{eq:objective}
    J(p;y^m) := \frac{1}{2}\| y(p)-y^m \|_{\Sigma_l^{-1}}^2+ \frac{1}{2}  \| p-p^0  \|^2_{\Sigma_p^{-1}} 
\end{equation}
over $p\in\mathcal{X}$, where $\left \| v \right \|_{\Sigma_{}^{-1}} := v^T \Sigma_{}^{-1} v$ for any given $\Sigma$.

Assuming that $y(p)$ is continuous, it can be shown that the entire functional $J$ is continuous. This, and compactness of the parameter space $\mathcal{X}$,  guarantee the existence of a minimizer. In general, the minimization problem~\eqref{eq:objective} is nonlinear and possibly non-convex, such that local minima corresponding to different, locally most likely, parameter sets can exist. We focus, however, on the setting of classical parameter identification, where the model is compatible with the data and, consequently, a good parameter choice approximates the true parameters well and leads to a small data mismatch $\|y(p(y^m))-y^m\|$ and locally unique minimizers.

Using gradient based methods such as the Gauss-Newton algorithm~\cite{Deuflhard2004} with steps $\Delta p(p,y^m)$ satisfying
\begin{equation}\label{eq:gauss-newton}
 A(p)\Delta p = y_p(p)^T\Sigma_l^{-1}(y(p)-y^m) + \Sigma_p^{-1}(p^0-p), \quad A(p) := y_p(p)^T \Sigma_l^{-1}y_p(p) + \Sigma_p^{-1},
\end{equation}
unconstrained local minima $p(y^m)\in\mathop\mathrm{int}\mathcal{X}$ of $J$ can be computed as long as $y_p\in\R^{m\times d}$ is of full rank in a neighborhood of the minimizer $p_{\rm MAP}$, which we assume throughout the paper. In the small residual setting assumed here, the Gauss-Newton algorithm exhibits a fast linear convergence.

The computation is efficient, if evaluations of $y(p)$ and the derivative $y_p(p)$ are replaced by approximations $y_\epsilon(p)$ and $y_{p;\epsilon}(p)$  with sufficiently small error $\epsilon$. Even if the number of forward problem evaluations necessary for computing a (local) maximum posterior point estimate $p_{\rm MAP}$ is vastly smaller than for posterior sampling approaches, the computational cost of evaluating several sufficiently accurate model responses $y_\epsilon(p)$ during the minimization can be prohibitively large in real-time applications such as quality control.

In such cases, the model $y$, or its numerical realizations $y_\epsilon$, can be replaced by a surrogate model approximating the original model sufficiently well. In the next section, we consider the use of Gaussian process regression for building such surrogate models.

\subsection{A Gaussian process surrogate model} \label{sec:surrogate}

We aim at building a surrogate model $y^*$ for $y:\mathcal{X} \to \R$ based on a set of model simulations with certain accuracies.

\paragraph{Designs and training data.} Let $D:=\{\mathcal{D}:\mathcal{X} \to \R_+\cup \{\infty\} \mid \card\dom\mathcal{D}\in\N\}$ denote the set of admissible designs, i.e. positive functions on $\mathcal{X}$ with an effective domain of finite cardinality. A design $\mathcal{D}$ comprises a set $P_t:=\dom\mathcal{D}$ of $n$ evaluation points $p_i\in\mathcal{X}$ with associated simulation errors $\epsilon_i = \mathcal{D}(p_i)<\infty$. 
Different designs can be compared in terms of the order relations $\le$ and $\ge$. Let $D'$ be given as a design in addition to $\mathcal{D}$, then $\mathcal{D}'\le\mathcal{D}$ means that the uncertainties $\epsilon_i'$ are less than or equal to $\epsilon_i$ for all $p_i$.

In combination with corresponding model evaluations $\hat y_i = y_{\epsilon_i}(p_i)$, a design $\mathcal{D}$ forms a complete training data set $\mathcal{D}_t :=  \left(\mathcal{D},(\hat y_i)_{i=1,\dots,n}\right)$. We will generally omit the subscript $t$ if the model evaluations are implicit from the context.

For brevity of notation, we focus on scalar model outputs $y\in\R$, i.e. $m=1$, since the vectorial case can be treated component-wise if measurement components are independent. 

In contrast to the bulk of literature~\cite{Chen,Forrester,Sacks}, we do not take the simulation results as ground truth, but include the evaluation tolerances $\epsilon_i$ explicitly into the training data. As a crude but deliberately simple model, we assume the actual evaluation errors $e_i := y_{\epsilon_i}(p_i)-y(p_i) \sim \mathcal{N}(0,\epsilon_i^2)$ to be independently and normally distributed.

\paragraph{Gaussian processes.}
Gaussian process regression is a powerful and versatile stochastic tool for function approximation~\cite{Candela}. A Gaussian process $\{X_p\}_{p\in\mathcal{X}}$ is a  collection of random variables such that the joint distribution of every finite subset of random variables is again a  multivariate Gaussian. The Gaussian process $\mathcal{GP}$ is completely defined by the mean $\mu(p)$ and and the pairwise covariance $k(p,p')$ for $p,p'\in\mathcal{X}$. 

Using a Gaussian process model we \textit{a priori} assume that the deterministic response $y(p)$ describing the underlying functional behaviour is a realisation of a random variable sampled from a Gaussian process
\begin{equation}\label{eq:GP}
    y\sim \mathcal{GP}(\mu,k).
\end{equation}
By considering $y-\mu$ instead of $y$ we may without loss of generality assume $\mu = 0$.

The kernel function $k$ contains information about the shape and structure, and in particular the smoothness, we expect the model $y$ to have~\cite{Duvenaud}.

We assume a stationary process, i.e. the kernel function to be translation-invariant: $k(p+q,p'+q)=k(p,p')$ for all $q\in \R^d$. This is somewhat restrictive, but simplifies the hyperparameter optimization or makes it tractable in the first place. Since we assume the model $y$ to be continuously differentiable, a covariance ensuring a reasonable degree of smoothness should be used, such as the Matérn kernel~\cite{Rasmussen} with order $\nu > 1$ or, simpler, the ubiquitous squared exponential 
\begin{equation}\label{eq:kernel}
k(p,p') = \sigma_f^2\exp\left(-\frac{1}{2}(p-p')^TL^{-1}(p-p')\right),
\end{equation}
which we will use here. The kernel function $k$ usually depends on a few hyperparameters that crucially determine the properties of the Gaussian process and therefore have to be chosen appropriately. In the case of the exponential kernel~\eqref{eq:kernel}, the hyperparameters comprise the spatial metric given by a symmetric positive definite matrix $L\in\R^{d\times d}$ and a scale factor $\sigma_f^2 \in\R$.

\paragraph{Inference.}
With the prior assumption of $y$ following a Gaussian process, and training data $\mathcal{D}$ given, a posterior probability density for the model outputs $y_i = y(p_i)$ for $1\le i\le n+1$ at the training parameter positions $p_i$ for $1\le i\le n$ and some inference parameter position $p_{n+1}\in\mathcal{X}$ can be defined in the Bayesian context as the product of prior and likelihood~\cite{Kuss}. The prior is derived from the Gaussian process assumption as
\[
    \pi_{\rm prior}(y) \propto \exp\left(-\frac{1}{2}y^TK^{-1}y\right)
\]
with the the symmetric positive definite covariance matrix $K_{ij} = k(p_i,p_j)$ for $i,j=1,\dots,n+1$. Assuming independent and normally distributed evaluation errors $e_i$, the likelihood is 
\[
    \pi_{\rm like}(y) \propto \exp\left(-\frac{1}{2}\sum_{i=1}^n (y_i-\hat y_i)^2 \epsilon_i^{-2}\right).
\]
Note that the likelihood does not depend on the inference value $y_{i+1}$. The posterior distribution is then given as
\[
    \pi_{\rm posterior}(y) = \pi(p\mid \mathcal{D}) \propto \pi_{\rm like}(y)\pi_{\rm prior}(y)
    \propto \exp\left(-\frac{1}{2} (y-\bar y)^T\Gamma^{-1} (y-\bar y) \right), 
\]
with posterior covariance $\Gamma = (K^{-1}+E^{-1})^{-1}$, $E^{-1}=\mathrm{diag}([\epsilon_1^{-2},\dots,\epsilon_n^{-2}, 0])$, and posterior mean 
\begin{equation}\label{eq:GPR-posterior-mean}
    \bar y=\Gamma E^{-1} [\hat y, 0]^T.
\end{equation}
Consequently, the conditional mean of the inference value is $\bar y_{n+1}$ with variance 
\begin{equation}\label{eq:variance-estimate}
    \sigma^2_{\mathcal{D}}(p_{n+1}) := \sigma_{n+1,n+1}^2 = \Gamma_{n+1,n+1}.
\end{equation}
The value $\bar y_{n+1}$ is the \emph{best linear unbiased predictor} (BLUP)~\cite{Park} of~\eqref{eq:GP} conditioned on the training data $\mathcal{D}$. Note that by precomputation of quantities depending only on the training data $\mathcal{D}$, an efficient evaluation of $\bar y_{n+1}$ and its derivative with respect to the inference position $p_{n+1}$ is possible~\cite{ShiChoi}. The thus defined mapping $p_{n+1}\mapsto \bar y_{n+1}$ will be denoted by $y_{\mathcal{D}}:\mathcal{X}\to\R$, and the standard deviation $p_{n+1}\mapsto \sigma_{n+1,n+1}$ with $\sigma_{\mathcal{D}}:\mathcal{X}\to\R$.

\begin{remark}
    The assumption of the mean $\mu$ be known a priori is quite strong in practice, and known as \emph{simple kriging}~\cite{Krige1951ASA}. More sophisticated GPR approaches make an appropriate ansatz for $\mu$ that is inferred from the data along with the values $y(p)$~\cite{UKRIGING}. We stick to simple kriging for ease of presentation.
\end{remark}

\paragraph{Estimating derivatives.}
Besides the possibility to predict the mean value $\bar y_{n+1}$, it is also possible to predict the derivative $\bar y_{p;n+1}:\mathbb{R}^d\rightarrow \mathbb{R}^d$ of the conditional mean at a given position p.

Let $\bar K = K_{ij}  \in \mathbb{R}^{n \times n}$ for $i,j=1,\dots,n+1$ the reduced covariance matrix. We also define $\bar E = \mathrm{diag}([\epsilon_1^{2},\dots,\epsilon_n^{2}]) \in \mathbb{R}^{n  \times n}$ as the reduced matrix of evaluation errors and $k(p,P_t)$ as the covariance between a sample point $p$ and the training data $P_t$. With the predictive mean $\bar y_{n+1} = k(p_{n+1},P_t) (\bar K+\bar E)^{-1} \hat y$, the derivative calculates as 
\begin{align}\label{eq:derivative}
    \frac{\partial \bar y_{n+1}(p)}{\partial p}\bigg|_{p=p_{n+1}} &= \frac{\partial k(p,P_t)}{\partial p}\bigg|_{p=p_{n+1}} \left(\bar K+\bar E\right)^{-1}\hat y.
\end{align}

It can be shown that the gradient follows a multivariate normal distribution, with the expected value corresponding to equation (\ref{eq:derivative}). This follows directly from the linearity of the derivative and the calculation of the expected value. 

\paragraph{Hyperparameter optimization.}
Hyperparamters $h \in S \subset (0,\infty)^{d+1}$ are free parameters within the kernel function determined from the training data $\mathcal{D}$. The choice of appropriate hyperparameters is crucial to ensure good predictive ability.
The hyperparameters are determined by minimising the negative log marginal likelihood \cite{gramacy2020surrogates,Rasmussen}, where the marginal likelihood (model evidence) is given by
\begin{align}
    \pi(\hat y|P_t,h) &= \int \pi_{\rm like}(\hat y|y)\pi_{\rm prior}(y)\;dy \nonumber \\
    &= \mathcal{N}(\hat y\,|\,0,\bar K+ \bar E).
\end{align}
with $\bar K = k(p_i,p_j)$ for $i,j=1,\dots,n$ and $\bar E=\mathrm{diag}([\epsilon_1^{2},\dots,\epsilon_n^{2}]) $.
The square exponential kernel (\ref{eq:kernel}) contains the hyperparameters $l = \mathrm{diag}([l_1,\dots,l_d]) \in\R^{d\times d} $ and $\sigma_f^2 \in\R$, where $l_i$ determines the horizontal dependence of the data in the corresponding feature direction and the scaling factor $\sigma_f$ determines the distance of the latent variables $f$ from the mean of the GP prior.
We therefore define the set $h$ consisting of $h := \left \{ \sigma_f^2,l_1,\dots,l_d \right \}$.
\begin{remark}
    The metric $l$ is chosen to be diagonal because of the simplicity and the small number of parameters, and since it allows an automatic relevance determination. On the other hand, the choice implies independence of parameter dimensions, which may be an incorrect model assumption and therefore lead to incorrect model results. For more complex models and their implications we refer to~\cite{Poggio}.
\end{remark}
There is no guarantee that the log marginal likelihood does not contain multiple local minima, since it is not convex. To obtain reasonable parameters for the underlying parameter space, we further restrict the minimization problem by using box constraints for all parameters.
\begin{align}
    \min_{ h\in S} \quad & -\log{\pi( \hat y | P_t,h)}  = \frac{1}{2}\hat y^T \left ( \bar K(h)+E \right )^{-1}\hat y+ \frac{1}{2}\log\left ( \det{ \left ( \bar K(h)+E \right )}\right )+ \frac{n}{2} \log\left ( 2\pi \right )\\
    \textrm{s.t} \quad & \sigma \in\left [  \sigma_{lb}, \sigma_{ub} \right ],l_1\in\left [  l_{1,lb}, l_{1,ub} \right ],\dots, l_d \in\left [ l_{d,lb}, l_{d,ub} \right ].  \nonumber  
\end{align}
Using gradient based optimization we need to calculate the derivatives as follows
\begin{align}
    \frac{\partial}{\partial h_i} \log\left ( p(\hat y | P_t,h) \right ) &= \frac{1}{2}\hat y^T \bar K^{-1} \frac{\partial \bar K}{\partial h_i}  \bar K^{-1} \hat y - \frac{1}{2}\mathrm{tr}\left ( \bar K^{-1} \frac{\partial \bar K}{\partial h_i}  \right ) \nonumber\\
    &= \dfrac{1}{2} \left ( (\alpha\alpha^T-\bar K^{-1})\dfrac{\partial \bar K}{\partial h_i}\right ) ,\;\; \alpha = \bar K^{-1}\hat y ,\; i = 0,\dots,d.
\end{align}

\section{Adaptive Gaussian Process Regression}\label{sec:active-learning}

Replacing exact model evaluations $y(p)$ in the objective~\eqref{eq:objective} by a cheaper surrogate model $y_{\mathcal{D}}(p)$ yields maximum posterior point estimates 
\[
p_{\mathcal{D}}(y^m) = \argmin_{p\in \mathcal{X}} J_{\mathcal{D}}(p;y^m)
:= \frac{1}{2}\|y_{\mathcal{D}}(p)-y^m\|_{\Sigma_l^{-1}}^2 + \frac{1}{2}\|p-p^0\|_{\Sigma_p^{-1}}^2
\]
and saves computational effort for computing Gauss-Newton steps $\Delta p_{\mathcal{D}}(p,y^m)$ by solving
\begin{equation*}
    \tilde{A}_{\mathcal{D}}(p) \Delta p_{\mathcal{D}} 
    = y_{\mathcal{D},p}^T\Sigma_l^{-1}(y_{\mathcal{D}}(p)-y^m) + \Sigma_p^{-1}(p^0-p), \quad \tilde{A}_{\mathcal{D}}(p) := y_{\mathcal{D},p}(p)^T\Sigma_l^{-1}y_{\mathcal{D},p}(p) + \Sigma_p^{-1}.
\end{equation*}
It also incurs both some error $p_{\mathcal{D}}(y^m)-p(y^m)$ of the resulting identified parameters and a considerable computational effort for evaluating the model according to $\mathcal{D}$ beforehand.

\begin{remark}
    In defining $p_{\mathcal{D}}(y^m)$ we assume that $J_{\mathcal{D}}(p;y^m)$ has a globally unique minimum. By restricting $p$ to a suitable neighborhood of a locally unique minimizer, all considerations here can be extended directly to more general settings.
\end{remark}

When unlimited computational resources are available, arbitrarily accurate simulations can be run to generate huge amounts of training data and achieve any desired accuracy. However, with a finite computational budget, the question immediately arises for which parameters $p_i$ simulations should be performed with which uncertainties $\epsilon_i$ to achieve the best accuracy.
This is a classical design of experiments problem for $\mathcal{D}$, with competing objectives of minimizing the expected surrogate model's approximation error $E(\mathcal{D})$ and minimizing the computational effort $W(\mathcal{D})$ for creating the training data. 
We consider the formulation
\begin{equation}\label{eq:doe-basic-problem}
    \min_{\mathcal{D}\in D} E(\mathcal{D}) \quad\text{subject to } W(\mathcal{D}) \le W_{\rm max},
\end{equation}
but could equivalently consider minimizing the total work while requesting a certain accuracy, or minimizing an arbitrary strict convex combination of work and error. In the following, we establish quantitative error estimates $E(\mathcal{D})$ and work models $W(\mathcal{D})$, and design a sequential greedy heuristic for solving~\eqref{eq:doe-basic-problem}.

\subsection{Accuracy model}

First we need to quantify the parameter reconstruction error $p_{\mathcal{D}}(y^m)-p(y^m)$ in terms of the measurement error variance $\Sigma_l$ and the surrogate model approximation quality $y_{\mathcal{D}}-y$ depending on the design $\mathcal{D}$.

\paragraph{Pointwise error estimates.}

We start by establishing an estimate of the parameter reconstruction error for deterministic functions $y(p)$ and $y_{\mathcal{D}}(p)$. The small residual assumption from Sec.~\ref{sec:inverse-problem} guarantees the local uniqueness of both the exact model reconstruction and the surrogate model reconstruction. As a by-product, this will also yield the unavoidable error level due to measurement noise.

\begin{theorem} \label{th:local-errorbound}
    Assume there are constants $0<\bar R,C_1,C_2<\infty$ and a parameter point $p^*\in\mathcal{X}$ such that the forward model $y$ satisfies the following conditions.
    \begin{enumerate}
        \item $y:B(p^*,\bar R) \to \R^m$ is twice continuously differentiable with bounded derivatives $\|y_p(p)\| \le C_1$  and $\|y_{pp}(p)\| \le C_2$ for all $p\in B(p^*,\bar R)$.
        Here, $B(p,r)$ denotes the open ball of radius $r$ around $p$.
        \item For some $y^m\in\R^m$, $p^* = p(y^m)\in\R^d$ is a minimizer of the objective $J(p,y^m)$  with small residual, i.e. 
        \begin{equation}\label{eq:small-residual}
            \left\| \Sigma_l^{-1}(y(p^*)-y^m)\right\| \le \frac{L}{3C_2} 
            \quad\text{with }
            L:=\lambda_{\rm min}\left(y_p(p^*)^T\Sigma_l^{-1}y_p(p^*)+\Sigma_p^{-1}\right) > 0,
        \end{equation}
        where $\lambda_{\rm min}(\cdot)$ denotes the minimal eigenvalue of a given matrix.
    \end{enumerate}

    Then, there are $\bar\epsilon>0$ and $0<\bar\epsilon' <L/(3\|\Sigma_l^{-1}\|C_1)$, such that for all $\epsilon\le \bar\epsilon$ and $\epsilon' \le \bar\epsilon'$ the bound 
    \begin{equation}\label{eq:radius-definition}
        R := \frac{3\|\Sigma_l^{-1}\|(\epsilon'+C_1)\epsilon + L \epsilon'/C_2}{L-3\|\Sigma_l^{-1}\|C_1\epsilon'} < \bar R
    \end{equation}
    holds and for all surrogate models $y_\mathcal{D}:\R^d\to\R^m$ with
    \begin{equation}  \label{eq:surrogate-model-deviation}
    \|y_{\mathcal{D}}-y\|_{L^\infty(B(p^*,R))} \le \epsilon  \quad\text{ and }\quad \|(y_{\mathcal{D}})_p - y_p\|_{L^\infty(B(p^*,R))} \le \epsilon' ,
    \end{equation}
    there is a locally unique minimizer $p_\mathcal{D}(y^m)$ of $J_\mathcal{D}$ satisfying the error bound
    \begin{equation}\label{eq:error-bound}
        \left\|p_{\mathcal{D}}(y^m)-p^*\right\| \le R.
    \end{equation}
\end{theorem}
\begin{proof}
    The bound~\eqref{eq:radius-definition} is trivially satisfied for sufficiently small $\bar\epsilon$ and $\bar\epsilon'$. 
    
    The proof of~\eqref{eq:error-bound} will be based on an implicit path $\psi(t)$ connecting $p^*$ and some $p_\mathcal{D}(y^m)$. For that, we first define a linear interpolation between the exact model $y$ and the surrogate model $y_\mathcal{D}$ as
    \[
        y^t(p) := (1-t)y(p)+t y_{\mathcal D}(p) \quad\text{for $t\in[0,1]$},
    \]
    as well as the normal equations
    \[
        F(t,p) := y_p^t(p)^T\Sigma_l^{-1}(y^t(p)-y^m)+\Sigma_p^{-1}(p-p^0)=0.
    \]
    The implicit function theorem yields a path $\psi(t)$ satisfying $F(t,\psi(t))=0$ in a neighborhood of $\psi(0)=p^*$ if $F_p(0,\psi(0))$ is invertible. The path satisfies the ordinary differential equation $\psi_t = F_p(t,\psi(t))^{-1} F_t(t,\psi(t))$. Its existence up to $t=1$ and the bound $\|\psi(1)-\psi(0)\| \le R$, which yields the claim~\eqref{eq:error-bound}, is guaranteed by the Theorem of Picard-Lindelöf~\cite{Aulbach} if 
    \begin{equation}\label{eq:Picard-Lindeloef-bound}
        \|F_p(t,p)^{-1}F_t(t,p)\| \le R
    \end{equation}
    holds for all $p\in \overline{B(p^*,R)}$ and $t\in[0,1]$. In the following, we will establish invertibility of $F_p$ and the bound~\eqref{eq:Picard-Lindeloef-bound} for sufficiently small $\bar\epsilon,\bar\epsilon'$.

    For all $p\in \overline{B(p^*,R)}$ we have
    \begin{align}
    F_p(t,p)
    &= (y_{pp}^t)^T\Sigma_l^{-1}(y^t-y^m) + (y_{p}^t)^T\Sigma_l^{-1}y^t_p + \Sigma_p^{-1} \notag \\
    &= (y_{pp}^t)^T\Sigma_l^{-1}(y^t-y^m) + t^2 (y_\mathcal{D}-y)_p^T \Sigma_l^{-1}(y_\mathcal{D}-y)_p + 2t (y_\mathcal{D}-y)_p^T \Sigma_l^{-1}y_p + y_p^T \Sigma_l^{-1}y_p + \Sigma_p^{-1} \notag \\
    &= (y_{pp}^t(p))^T \left(\Sigma_l^{-1}(y^t(p)-y(p)+y(p)-y(p^*)+y(p^*)-y^m)\right) \label{eq:est-1} \\
    &\quad + t^2 (y_\mathcal{D}-y)_p^T \Sigma_l^{-1}(y_\mathcal{D}-y)_p + 2t (y_\mathcal{D}-y)_p^T \Sigma_l^{-1}y_p \label{eq:est-2} \\
    &\quad + (y_p(p)-y_p(p^*))^T\Sigma_l^{-1}(y_p(p)-y_p(p^*))
           + 2(y_p(p)-y_p(p^*))^T\Sigma_l^{-1}y_p(p^*) \label{eq:est-3}\\
    &\quad + y_p(p^*)^T\Sigma_l^{-1}y_p(p^*) + \Sigma_p^{-1}, \notag
    \end{align} 
    and can bound the norms of the individual terms as
    \begin{align*}
        \|\eqref{eq:est-1}\| 
        &\le C_2 \|\Sigma_l^{-1}(y^t(p)-y(p)+y(p)-y(p)+y(p)-y^m)\| \\
        &\le C_2\left ( \left \| \Sigma_l^{-1}\left (  y^t(p)-y(p)\right ) \right \| + \left \| \Sigma_l^{-1}\left (  y(p)-y(p)\right )      \right \| + \left \| \Sigma_l^{-1}\left (  y(p)-y^m\right ) \right \| \right )  \\        
        &\le C_2 \|\Sigma_l^{-1}\|(\epsilon + RC_1) + \frac{L}{3}, \\[1ex]
        \|\eqref{eq:est-2} \|
        &\le \|\Sigma_l^{-1}\| \left( (\epsilon')^2  + 2\epsilon'C_1\right), 
        \quad\text{and}\\[1ex]
        \|\eqref{eq:est-3}\|
        &\le \|\Sigma_l^{-1}\| (R^2 + 2RC_1)
    \end{align*}
    by using Taylor's theorem, the small residual assumption~\eqref{eq:small-residual}, and $t\le 1$. 
    Consequently, we obtain
    \[
        \lambda_{\rm min}(F_p(t,p)) \ge \frac{2}{3}L
        - \|\Sigma_l^{-1}\|\left( C_2\epsilon + (C_2+2)RC_1 + (\epsilon')^2  + 2\epsilon'C_1 + R^2  \right).
    \]
    Since $R = \mathcal{O}(\epsilon + \epsilon')$, there exist sufficiently small $\bar\epsilon,\bar\epsilon'>0$, such that for all $\epsilon<\bar\epsilon,\epsilon'<\bar\epsilon'$
    \[
        \lambda_{\rm min}(F_p(t,p)) \ge \frac{L}{3}
    \]
    holds. Consequently, $F_p$ is positive definite with $\|F_p(t,p)^{-1}\| \le 1/ \lambda_{\rm min} \le 3/L$.
    
    Moreover, we obtain
    \begin{align*}
    \|F_t(t,p)\|
    &= \| (y_\mathcal{D}-y)_p^T \Sigma_l^{-1}(y^t-y^m) + (y^t_p)^T \Sigma_l^{-1}(y_\mathcal{D}-y) \| \\
    &=\| (y_\mathcal{D}-y)_p^T \Sigma_l^{-1}(y^t-y+y-y(p)+y(p)-y^m) + (y^t_p)^T \Sigma_l^{-1}(y_\mathcal{D}-y) \| \\
    &\le \epsilon' \left(\|\Sigma_l^{-1}\|(\epsilon+RC_1) + \frac{L}{3C_2}\right) + C_1 \|\Sigma_l^{-1}\| \epsilon\\
    &\le \|\Sigma_l^{-1}\|(\epsilon'(\epsilon + RC_1)+\epsilon C_1) + \frac{\epsilon' L}{3C_2}
    \end{align*}
    and thus
    \begin{align}
    \| F_p^{-1} F_t\| 
    &\le \frac{3\|\Sigma_l^{-1}\|}{L} (\epsilon'(\epsilon + RC_1)+\epsilon C_1) + \frac{\epsilon'}{C_2}.  \label{eq:righthandside}
    \end{align}
    Finally, we bound the right hand side of~\eqref{eq:righthandside} by $R$ using~\eqref{eq:radius-definition} via
    \begin{align*}
        & R = \frac{3\|\Sigma_l^{-1}\|(\epsilon'+C_1)\epsilon + L \epsilon'/C_2}{L-3\|\Sigma_l^{-1}\|C_1\epsilon'} \\
        \Rightarrow \quad & 3\|\Sigma_l^{-1}\| RC_1 \epsilon' -RL = -\frac{\epsilon' L}{C_2} - 3\|\Sigma_l^{-1}\|\epsilon C_1 - 3\|\Sigma_l^{-1}\|\epsilon \epsilon' \\
        \Rightarrow \quad & 3\|\Sigma_l^{-1}\| (\epsilon'(\epsilon + RC_1)+\epsilon C_1) + \frac{\epsilon'}{C_2} L = RL\\
        \Rightarrow \quad & \frac{3\|\Sigma_l^{-1}\|}{L} (\epsilon'(\epsilon + RC_1)+\epsilon C_1) + \frac{\epsilon'}{C_2} = R,
    \end{align*}
    and obtain $\|F_p^{-1}F_t\| \le R$, i.e. the required estimate~\eqref{eq:Picard-Lindeloef-bound}.

    By construction, $\psi(1)\in\overline{B(p^*,R)}$ is a stationary point of $J_{\mathcal{D}}$ due to $(J_{\mathcal{D}})_p(\psi(1);y^m) = F(1,\psi(1)) = 0$. Since $F_p(1,p)=(J_{\mathcal{D}})_{pp}$ is positive definite for all $p\in\overline{B(p^*,R)}$, we conclude that $p_{\mathcal{D}}(y^m) := \psi(1)$ is a locally unique local minimizer of $J_{\mathcal{D}}$ satisfying~\eqref{eq:error-bound}.
\end{proof}

\begin{remark}
    Equation \eqref{eq:surrogate-model-deviation} takes into account not only the approximation of the model $y$ by the surrogate model but also the adequate approximation of the derivatives of the model $y'$. It is assumed, and numerical experiments support this assumption, that if the error of the model approximation is small, the error of the derivatives is also small.
\end{remark}

The error bound~\eqref{eq:radius-definition} can be simplified without loosing asymptotic accuracy regarding $\epsilon, \epsilon' \rightarrow 0$, clearly revealing the linear dependence of the parameter reconstruction error on the surrogate model accuracy.

\begin{corollary}
    Let the assumptions of Theorem~\ref{th:local-errorbound} be satisfied. Then there are $\bar\epsilon,\bar\epsilon'$ such that the claim~\eqref{eq:error-bound} also holds for 
    \begin{equation}\label{eq:linear-error-bound}
        R = 12L^{-1}\|\Sigma_l^{-1}\|C_1\epsilon + C_2^{-1}\epsilon'.
    \end{equation}
\end{corollary}
\begin{proof}
    By restricting $\bar\epsilon'\le \min\{L/(3\|\Sigma_l^{-1}\|C_1),C_1\}$ we obtain $L-3\|\Sigma_l^{-1}\|C_1\epsilon' \le 1/2$ and, inserting this into~\eqref{eq:radius-definition}, directly proves the claim.
\end{proof}

We point out, that Thm.~\ref{th:local-errorbound}, while establishing the stability structure present in the approximate parameter identification problem, can usually not be applied directly for numerical computation, since its assumptions are hard to verify in practice. For the construction of actual algorithms below, we will therefore rely on computable estimates following the structure provided by the above theory.

\paragraph{Relevant error quantity.}

We assume that a small absolute parameter reconstruction error is desired, but, since there is some unavoidable error due to measurement errors, a small relative error is also sufficient. First we estimate the unavoidable error level $e_0$.

\begin{corollary}
    Let the assumptions of Theorem~\ref{th:local-errorbound} be satisfied. Let $p(y^m)$ be a locally unique minimizer of $J(p,y^m)$ and $\delta \in \R^m$ some measurement noise. Then there is some $\bar\epsilon>0$, such that for $\|\delta\| \le \bar\epsilon$ there is a locally unique minimizer $p(y^m+\delta)$ of $J(p,y^m+\delta)$ with 
    \[
        \|p(y^m)-p(y^m+\delta)\| \le \frac{3C_1\|\Sigma_l^{-1}\|}{L} \|\delta\|.
    \]
\end{corollary}
\begin{proof}
    We define the auxiliary surrogate model $y_\mathcal{D} := y-\delta$, which obviously satisfies~\eqref{eq:surrogate-model-deviation} with $\epsilon = \|\delta\|$ and $\epsilon' = 0$. Note that the corresponding objective $J_\mathcal{D}$ satisfies $J_\mathcal{D}(p,y^m) = J(p,y^m+\delta)$. Theorem~\ref{th:local-errorbound} guarantees the existence of a local minimizer $p_\mathcal{D}(y^m) = p(y^m+\delta)$ satisfying~\eqref{eq:error-bound}, which proves the claim.
\end{proof}

With an expected measurement error magnitude $\|\delta\|$ of $\sqrt{\|\Sigma_l\|}$, we define the unavoidable error level as 
\[ 
  e_0 := \frac{3C_1}{L} \|\Sigma_l^{-1}\| \sqrt{\|\Sigma_l\|}.
\]
Since $L$, $\epsilon$, $\epsilon'$, and, in principle, also $C_1$ and $C_2$ depend on $p=p(y^m)$, and $\epsilon,\epsilon'$ also on the surrogate model specified by the design $\mathcal{D}$, the values of $e_0$ and $R$ can be localized, and we write $R_\mathcal{D}(p)$ and $e_0(p)$ explicitly.

Aiming at a low absolute error while allowing for a certain relative error, we define the local error quantity 
\begin{equation}\label{eq:local-error-model}
e_\mathcal{D}(p) := \frac{R_\mathcal{D}(p)}{1+\alpha e_0(p)} 
                 \le \max\left\{ R_\mathcal{D}(p), \alpha \frac{R_\mathcal{D}(p)}{e_0(p)}\right\}
\end{equation}
that is to be minimized by selecting an appropriate design $\mathcal{D}$.
Here, $\alpha\ge 0$ acts as an arbitrary weighting factor of absolute and relative accuracy.

Since during the construction of the surrogate model $y_\mathcal{D}$ by minimizing~\eqref{eq:doe-basic-problem} the measurement values $y^m$ and hence the parameter position $p=p(y^m)$ of interest are unknown, the error quantity~\eqref{eq:local-error-model} needs to be considered over the whole parameter region $\mathcal{X}$.
We therefore define the accuracy model 
\begin{equation}\label{eq:error-model}
    \mathcal{E}(\mathcal{D}) := \|e_\mathcal{D}\|_{L^q(\mathcal{X})} \quad \text{for some $1\le q < \infty$}.
\end{equation}
Choosing $q\approx 1$ would focus on minimizing the average parameter reconstruction error, while choosing $q$ very large would focus on the worst case. Note that for $q<\infty$ the accuracy model is continuously differentiable in $\epsilon,\epsilon'$ as functions on $\mathcal{X}$. For the numerical experiments in Sec.~\ref{sec:experiments} we have chosen $q=2$.

Still missing are the surrogate model error bounds $\epsilon,\epsilon'$ in terms of the design $\mathcal{D}$. Unfortunately, for virtually all cases of practical interest, there is little hope for obtaining simultaneously rigorous and quantitatively useful bounds. For GPR surrogate models in particular, the global support of the posterior probability density precludes the existence of a strict bound, though its fast decay provides thresholds that are not exceeded with high probability. The assumed normal distribution of errors $y(p)-y_\mathcal{D}(p) \sim \mathcal{N}(0,\sigma_{\mathcal{D}}(p)^2)$ implies that $\|y(p)-y_\mathcal{D}(p)\|$ is generalized-$\chi^2$-distributed and hence formally unbounded. Instead of a strict bound, one can use a representative statistical quantity for $\epsilon$, such as the median $\epsilon := \sqrt{m(1-2/(9m))^3}\|\sigma_\mathcal{D}(p)\|_2$ of a corresponding $\chi$ distribution as an upper bound for the square root of the generalized $\chi^2$ distribution's median, or the mean $\epsilon^q := \mathop\mathrm{tr}(\sigma_\mathcal{D}(p)^q)$ of the generalized $\chi^2$ distribution itself~\cite{MathaiProvost1992}.

Numerical experience suggests a roughly proportional behaviour of $\epsilon$ and $\epsilon'$, with a proportionality constant $\beta$ that can be estimated from available model evaluations. We therefore assume $\epsilon'=\beta\epsilon$ and simplify the error bound~\eqref{eq:linear-error-bound} accordingly. 

Inserting the thus chosen values of $\epsilon$ and $\epsilon'$ into~\eqref{eq:linear-error-bound} completes the accuracy model. All together, the error model assumes the form of a weighted $L^q$ norm with weights $w(p)$ that can be explicitly computed, either directly or as an estimate:
\begin{equation}\label{eq:globalerror_integral}
    \begin{aligned}
    E(\mathcal{D}) &= \left(\int_{p\in\mathcal{X}} w(p)^q \epsilon(p)^q \, dp\right)^{\frac{1}{q}} \approx \mathcal{E}(\mathcal{D}),  \\
    w(p)^q &= \left ( \left ( 1+\alpha e_0(p) \right )^{-1}\left ( 12L^{-1}\left \| \Sigma_l^{-1} \right \|C_1+\beta C_2^{-1} \right ) \right )^q, \\
    \epsilon(p)^q &= \begin{cases} 
        (m(1-\frac{2}{9m})^3)^{\frac{q}{2}} \|\sigma_\mathcal{D}(p)\|_2^q, &\text{or} \\
        \mathop\mathrm{tr}(\sigma_\mathcal{D}(p)^q),
        &\text{depending on the choice above}.
        \end{cases}
    \end{aligned}
\end{equation}

\paragraph{Numerical integration.}
The integral in (\ref{eq:globalerror_integral}) can be computed analytically only in the rarest cases, so we can rely on the approximation by Monte Carlo integration. This approximation is easy to implement for higher-dimensional parameter spaces and greatly simplifies the subsequent computation of the gradient and Hessian for the later optimization. The estimate is given by
\begin{equation}\label{eq:MC-integration}
  E(\mathcal{D}) \approx  \left ( \frac{\mathrm{vol}(\mathcal{X})}{{N_{\rm MC}-1}}\sum_{i=0}^{N_{\rm MC}-1}w^q(p_i) \epsilon(p_i)^q  \right )^{\frac{1}{q}}
\end{equation}
where $N_{\rm MC}$ is the number of sample points used to evaluate the above expression.
\begin{remark}
    To obtain a reliable value for the integral approximation, a sufficient number of points $N_{\rm MC} > 0 $ in the parameter space $\mathcal{X}$ must be evaluated \cite{Gamerman}. The large number of evaluation points is not problematic, since the weights $w^q(p_i)$ as well as the standard deviation $\sigma_\mathcal{D}(p_i)$ can be calculated with little computational effort.
\end{remark}
\paragraph{Weight factors.}
Since it is most likely impossible to obtain a priori information about $C_1$ and $C_2$ used in Thm.~\ref{th:local-errorbound}, we cannot explicitly compute the weights $w^q(p_i)$. Instead, we can compute estimates $\tilde w$ for them directly from the Gauss-Newton method applied to $F(p)=\left \|f(p)\right \|^2_{\Sigma^{-1}_l} =\left \| y(p)-y_{\mathcal{D}}(p) \right \|^2_{\Sigma^{-1}_l}\rightarrow \min $. Performing a Gauss Newton step yields
\begin{equation}\label{eq:parametererrorestimate}
    p_{\mathcal{D}}-p  \approx \left ( f'^T \Sigma_l^{-1} f' \right )^{-1}f'^T\Sigma_l^{-1} \left( y_{\mathcal{D}}-y \right)
\end{equation}
with the Jacobian $f' = f'(p)\in \R^{m \times d}$ and thus the weight $\tilde{w}$ as the error transport factor
\begin{equation}\label{eq:weightfactordefintion}
    \tilde{w}(p) := \left \| (f'^T \Sigma_l^{-1}f'+\lambda I)^{-1}f'^T\Sigma_l^{-1} \right \|_2.
\end{equation}
Regularisation with $\lambda \in \mathbb{R}_{>0}$ may be necessary if the row-wise entries of the Jacobian are very similar, which is the case when the experimental data on which the m regression models are based are similar.

Note that the derivatives $f'(p)$ correspond to the estimated derivatives according to \eqref{eq:derivative} with or without gradient data corresponding to the current data set.

\subsection{Work model}

The evaluation of the forward model $y(p)$ usually involves some kind of numerical approximation such as discretization of differential equations, iterative solution of equation systems, or Monte Carlo sampling of stochastic systems, which results in an approximation $y_\epsilon(p)$. While in principle any uncertainty $\|y_\epsilon(p)-y(p)\| \le \epsilon$ for arbitrary $\epsilon>0$ can be achieved, the accuracy requires a computational effort $W(\epsilon)$ to be spent on the evaluation. 

First we recall from~\cite{Weiser} two prototypical work models for the settings that are most relevant due to the large effort of the involved computations.

\paragraph{Finite element discretization.} Let us assume the value of $y$ is extracted from a solution of an elliptic partial differential equation on a domain in $\R^d$, and the solution is approximated by finite elements of order $r\ge 1$ on an adaptively refined mesh with $N$ vertices. Then, the discretization error $\epsilon$ can be expected to be proportional to $N^{-r/d}$, see~\cite{DeuflhardWeiser2012}. With a solver of optimal complexity, such as multigrid with nested iteration, the computational work for obtaining a solution of uncertainty $\epsilon$ is of order $N$, and we obtain
\begin{equation}\label{eq:fe-work}
    W(\epsilon) = \frac{r}{d}\epsilon^{-d/r}.
\end{equation}
The factor $r/d$ is introduced for later convenience, and does not affect the selection of evaluation positions or accuracies in any way.

If a sparse direct solver is used instead, i.e., for time-harmonic Maxwell equations, the computational work is rather on the order of $N^{1.5}$ for $d=3$~\cite{George}, which leads to 
\[
    W(\epsilon) = \frac{r}{1.5d}\epsilon^{-1.5d/r}.
\]
Of course, the asymptotic behavior $W\to 0$ for $\epsilon\to\infty$ is not realistic, as there is a fixed amount $W_{\rm min}$ of work necessary on the coarsest grid. Thus, the work model is strictly valid only for $\epsilon \le \epsilon_{\rm max}$.

\paragraph{Monte Carlo sampling.}
In case the forward model $y$ contains a high-dimensional integral to be evaluated by Monte-Carlo sampling, standard convergence results suggest a sampling error $\epsilon$ proportional to $N^{-1/2}$, where $N$ is the number of samples. This leads to the work model
\[
    W(\epsilon) = \frac{1}{2}\epsilon^{-2},
\]
which is just a special case of~\eqref{eq:fe-work}.

Omitting constant factors without loss of generality, we therefore assume that the computational work is given by 
\begin{equation}\label{eq:generic-work-model}
    W(\epsilon) = \epsilon^{-2s}, \quad s > 0
\end{equation}
in the remainder of the paper.

\begin{remark}
    The work models \eqref{eq:fe-work} and \eqref{eq:generic-work-model} are highly idealized and based on asymptotic theoretical results for $\epsilon\to 0$. The actual computational effort depending on a tolerance $\epsilon$ is instead piecewise constant and, in particular for large tolerances, rather different from the model. Fortunately, the models are reasonably accurate for small tolerances, where the most effort is spent, and therefore useful for design optimization.
\end{remark}

The computational effort for a complete training data set $\mathcal{D}$ based on the design $\mathcal{D}$ is just the sum of the individual simulation efforts, i.e.
\[
    W(\mathcal{D}) := \sum_{p\in\dom\mathcal{D}} W(\mathcal{D}(p)).
\]
Here, we simply neglect the highly problem-specific dependence of the computational effort on the evaluation parameter $p$. We note that $W(\mathcal{D})$ inherits convexity, monotonicity, and the barrier property on $\R^n_+$ from the individual contributions $W(\epsilon_i)$.

We will also be interested in the computational effort coming with \emph{incremental designs}. Assume $\mathcal{D}$ is a design that has already been realized. Evaluating the model on a finer design $\mathcal{D}' \le \mathcal{D}$ can consist of simulating the model for parameters $p\not\in\dom\mathcal{D}$, or improving the accuracy of already performed simulations for $p\in\dom\mathcal{D}$ with $\mathcal{D}'(p)< \mathcal{D}(p)$, or both. If already conducted simulations have been stored such that they can be continued instead of started again, the computational effort of obtaining the training data set $\mathcal{D}'$ from $\mathcal{D}$ is
\[
    W(\mathcal{D}'|\mathcal{D}) = W(\mathcal{D}') - W(\mathcal{D}).
\]

\subsection{The design of computer experiments problem}\label{ss:computer experiment}

In order to create a surrogate model $y_{\mathcal{D}}$ as discussed in Sec.~\ref{sec:surrogate}, an appropriate design $\mathcal{D}$ must be found by solving the design of experiments problem~\eqref{eq:doe-basic-problem}. Since little is known about the model derivative $y_p$, and consequently about $E(\mathcal{D})$, before any simulations have been performed,~\eqref{eq:doe-basic-problem} cannot be reasonably solved a priori. We therefore follow a sequential design of experiments approach by incrementally spending computational budget. In each step, we thus have to solve the modified problem for an incremental design $\mathcal{D}'\in D$ refining a given preliminary design $\mathcal{D}$:
\begin{equation}\label{eq:doe-incremental}
    \min_{\mathcal{D}'\le \mathcal{D}} E(\mathcal{D}') \quad\text{subject to } W(\mathcal{D}'|\mathcal{D}) \le \Delta W.
\end{equation}
The design problem~\eqref{eq:doe-incremental} is of combinatorial nature due to the unknown number $n$ of evaluation points and the choice between introducing new points or re-using existing ones, highly nonlinear due to the parameter locations $p_i$ being optimized, and therefore difficult to address rigorously. In particular, a relaxation of the design to the space of nonnegative regular Borel measures as used in~\cite{NeitzelPieperVexlerWalter2019} is infeasible due to the nonlinearity of the work model~\eqref{eq:fe-work}.

We therefore simplify the problem~\eqref{eq:doe-incremental} by separating the selection of evaluation positions $p_i$ from the choice of evaluation accuracies $\epsilon_i$.

\subsubsection{Choice of candidate points}

Instead of directly optimizing the essential support of the refined design $\Tilde{\mathcal{D}}'$, i.e. the number and position of evaluation points, we aim at selecting the best points from a fixed set used as a discrete approximation of $\mathcal{X}$. Due to the constraint $\Tilde{\mathcal{D}}'\le \Tilde{\mathcal{D}}$, the already existing evaluation points need to be retained, such that $\dom \Tilde{\mathcal{D}}' \supset \dom\Tilde{\mathcal{D}}$ can be defined in terms of a set $\mathcal{P}_C = \left \{ p_{c,i}\;| \; p_{c,i}\in\mathcal{X}\backslash\dom\Tilde{\mathcal{D}},\; i=1,\dots,k \right \}$ of additional candidate points. We then rely on the accuracy optimization step discussed below to produce a good -- often sparse -- solution with $\dom\Tilde{\mathcal{D}}' \subset \dom\Tilde{\mathcal{D}} \cup \mathcal{P}_C$, i.e. selecting some additional evaluation points to be actually included into $\Tilde{\mathcal{D}}'$. 

Different ways of choosing the candidate points $p_{c,i}$ are conceivable. We can simply take $k$ random points in $\mathcal{X}$, either uniformly distributed or sampled from a density reflecting the error propagation factor $w(p)$. Alternatively, candidate points can be taken from low discrepancy sequences such as the Halton sequence~\cite{Halton1960}. 
When selecting points randomly or by sequences such as Halton or Sobol, the set of points can be filtered in advance for those points that have little or no effect on the global error $E$, which speeds up the optimization. As a criterion, the evaluation of the gradient $\nabla_{\epsilon} E$ of the global error can be used. Points for which
\begin{equation}
    \left |\nabla_{\epsilon} E( \dom\mathcal{D}') \right |_i < \mathrm{TOL}
\end{equation}
is component wise valid are then removed from the set of candidate points.

With larger computational effort, particularly promising candidate points can be obtained by finding local maximizers of the local error estimate density
\begin{equation}\label{eq:utility}
    g(p):=\Tilde{w}(p)\epsilon(p) = \Tilde{w}(p)\sqrt{\mathop\mathrm{tr}(\sigma_{\Tilde{\mathcal{D}}}(p)^2)}
\end{equation}
in analogy to the utility functions from Bayesian optimization~\cite{Pourmohamad}.

If candidate points are selected that are already in the training data, the next worst point determined by the acquisition function is selected to avoid re-adding.
In any case, points that are too close to each other should also be excluded in order to avoid numerical instabilities in GP evaluation.

The number $k$ of candidate points should be chosen carefully. A severely limited choice of candidate points due to small $k$ will likely lead to the inclusion of sub optimal points into $\dom\Tilde{\mathcal{D}}'$, whereas a large $k$ incurs a high computational effort in optimizing the evaluation errors $\epsilon_{c,i}$.

\subsubsection{Optimizing evaluation accuracies}

With a fixed set of candidate points to consider, the optimization problem~\eqref{eq:doe-incremental} is reduced to a nonlinear programming problem for the evaluation uncertainties $\epsilon_i$. We will, however, reformulate it equivalently in terms of the $n+j$ auxiliary variables $v_i := \epsilon_i^{-2}$, since then the objective $E$ is convex in $v$, see Thm.~\ref{th:convex-objective} below. Thus, we obtain
\begin{equation}
  \begin{aligned}\label{eq:doe-incremental-reduced}
    \min_{v \in \mathbb{R}^{n+j}_{+}} \quad &  E(v) \\ 
     \mathrm{s.t.}\quad& W(v) \leq \Delta W + W(\mathcal{D}),\\ 
    &v \geq \underline{v},
  \end{aligned}
\end{equation}
with
\begin{equation}\label{eq:transformed-objective}
   E(v) = \left(\int_{p\in\mathcal{X}} w(p)^q \epsilon(p)^q \, dp\right)^{\frac{1}{q}}, \quad
    \underline{v} = \begin{bmatrix} \epsilon^{-2} \\ 0 \end{bmatrix},\quad
    W(v) = \sum_{i=1}^{n+j} v_i^s
\end{equation}
due to~\eqref{eq:generic-work-model}, and $\sigma_v(p)$ given by~\eqref{eq:variance-estimate}.
Note that due to the change from $\epsilon$ to $v$, the upper bound $\mathcal{D}'\le\mathcal{D}$ in~\eqref{eq:doe-incremental} is transformed into the lower bound $v\ge\underline{v}$ in~\eqref{eq:doe-incremental-reduced}. To simplify the optimization, we minimize $\tilde{E}(v):=E(v)^q$ instead of $E(v)$, and  use the Monte Carlo approximation~\eqref{eq:MC-integration} of the integral. This leads to the minimization of
\begin{equation}
    \tilde{E}(v) \approx   \frac{\mathrm{vol}(\mathcal{X})}{{n+j-1}}\sum_{i=0}^{n+j-1}\tilde{w}(p_i) \epsilon(p_i)^q   =  \frac{\mathrm{vol}(\mathcal{X})}{{n+j-1}}\sum_{i=0}^{n+j-1}\tilde{w}(p_i)\mathop\mathrm{tr}(\sigma_\mathcal{D}(p_i)^2).
\end{equation}
under the same constraints as in \eqref{eq:doe-incremental-reduced} with $\tilde{w}(p_i) := w^q(p_i)$.

\begin{theorem} \label{th:convex-objective}
    For an exponent $q\ge 2$ in the error model~\eqref{eq:transformed-objective}, the objective $E(v)$ is convex. If, in addition, the covariance kernel $k(\cdot,\cdot)$ is strictly positive, $E$ is strictly convex and for any minimizer of~\eqref{eq:doe-incremental-reduced}, the work constraint $W(v)\le \Delta W + W(\mathcal{D})$ is active.
\end{theorem}
\begin{proof}
    Since the composition $(f\circ g)(x)$ of a convex function $g:\R^{n+j}\to\R$ and a convex, monotonically increasing $f:\R\to\R$ is quasi-convex, we only need to prove that the variance $\sigma^2_v(\tilde p)$ is convex in $v$ for all $p$.
    
    Let $V = \textrm{diag}([v, 0])$ such that the posterior covariance at point $\tilde p$ is given by
    \[
     \sigma^2_v(\tilde p) = e^T_{n+j+1} \left[K(\tilde p)^{-1}+V \right]^{-1} e_{n+j+1}
    \] 
    with the symmetric positive semi definite prior covariance $K(p)$ and the euclidean unit vector $e_{n+j+1}$. 
    
    For a small perturbation $0\ne\delta v\in\R^{n+j}$ and a corresponding $\delta V = \mathrm{diag}([\delta v, 0])$, we obtain for $K=K(\tilde p)$
    \begin{align*}
    \left[K^{-1}+V+\delta V\right]^{-1}
    &= \left[\left ( K^{-1}+V  \right )\left (I+ \left ( K^{-1}+V  \right )^{-1} \delta V \right )\right]^{-1} \\
    &= \left (I+ \left ( K^{-1}+V  \right )^{-1}  \delta V \right )^{-1} \left ( K^{-1}+V  \right )^{-1}.
    \end{align*}
    Introducing $A:=\left ( K^{-1}+V  \right )^{-1}$, the Neumann series representation yields
    \begin{align*}
    \left[K^{-1}+V+\delta V\right]^{-1}
    &= \left (I+\delta VA \right )^{-1}A \\
    &= \left (I-\delta VA +(\delta V A)^2 \right )A +\mathcal{O}(\delta V^3) \\
    &= \left (I-\delta VA +\delta V A\delta V A \right )A +\mathcal{O}(\delta V^3) \\
    &= A - A\delta VA  +  A \delta V A\delta V A  + \mathcal{O}(\delta V^3). 
    \end{align*} 
    Consequently, the second directional derivative of $\sigma^2_v(\tilde p)$ in direction $\delta v$ is
    \begin{equation}\label{eq:second-directional-derivative}
        2e_{n+j+1}^T A\delta VA\delta VA e_{n+j+1} \ge 0.
    \end{equation}
    Since this holds for any $\delta v$, the Hessian of $\sigma_v^2$ is positive semidefinite everywhere, and $\sigma_v^2$ therefore convex.
    
    Next we will show that a strictly positive covariance kernel $k(\cdot,\cdot)$ implies a positive $A$, which is the basis for proving the remaining claims.
    
    With $k$ strictly positive, the prior covariance $K$ is positive, i.e $K_{ij}>0~\forall i,j$, and thus $K^{-1}$ is an M-matrix, see~\cite{Plemmons}, and can be written as 
    \[
        K^{-1} = sI - B \quad \text{for some } B\ge 0, \quad s\ge \|B\|_2.
    \]
    We now observe that for $r=\|V\|_2 > 0$,
    \begin{align*}
        A^{-1} 
        &= K^{-1} + V \\
        &= sI - B + V \\
        &= (s+r)I - (B+rI-V)
    \end{align*}
    holds, with $rI- V\ge 0$ and $s+r \ge \|B\|_2 + \|rI- V\|_2 \ge \|B+rI- V\|_2$. Consequently, $A^{-1}$ is an invertible M-matrix, such that $A$ is nonnegative. Moreover, by the Neumann series, we can write
    \begin{align*}
        A 
        &= \left( (s+r)I - (B+rI-V)\right)^{-1} \\
        &= (s+r)^{-1} \sum_{i=0}^\infty (s+r)^{-i} (B+rI-V)^i \\
        &\ge (s+r)^{-1} \sum_{i=0}^\infty (s+r)^{-i} B^i.
    \end{align*}
    Since $K = s^{-1} \sum_{i=0}^\infty s^{-i}B^i$ is positive, so is $A$.
    
    If $A$ is positive, then so is $Ae_{n+j+1}$, and consequently $\delta V Ae_{n+j+1} \ne 0$, since $\delta V$ has at least one non vanishing entry. Then, the second directional derivative of $\sigma_v^2(\tilde p)$ as given in~\eqref{eq:second-directional-derivative} is  positive and $\sigma_v^2(\tilde p)$ strictly convex in $v$.  
    
    Moreover, the first directional derivative in direction $e_i$ is 
    \[
        \frac{\partial \sigma_v^2(\tilde p)}{\partial v_i} 
        = - e_{n+j+1}^T A e_i e_i^T A e_{n+j+1} = -A_{i,n+j+1}^2 < 0,
    \]
    and thus $\nabla_v \sigma_v^2(p) < 0$. Therefore, the first order necessary conditions can only be satisfied if the work constraint is active.
\end{proof}

The convexity of the admissible set $V:=\{v\in\R^{n+j}\mid v\ge \underbar{v} \wedge W(v)\le \Delta W + W(\mathcal{D})\}$ depends on the exponent $s$ showing up in the generic work model~\eqref{eq:generic-work-model}. Clearly, for $s\ge 1$, $W$ is convex, whereas for $s<1$ it is in general non convex, not even quasi-convex. In combination with Thm.~\ref{th:convex-objective}, we obtain the following result.

\begin{corollary}
    For exponents $q\ge 2$ and $s\ge 1$, the tolerance design problem~\eqref{eq:doe-incremental-reduced} is convex. If, in addition, the covariance kernel $k(\cdot,\cdot)$ is strictly positive, the design problem has a unique solution.
\end{corollary}

Remember that $q$ defines the type of norm in which the global error bound is measured. If $q=1$, it refers to the expected error over all possible points in $\mathcal{X}$, whereas $q=2$ implies a mean squared error notion, and $q\to\infty$ considers the maximum error. The choice of $q$ is therefore a modeling question of what error distribution is acceptable. For the sake of simplicity of exposition, we will restrict the attention to $q\ge 2$, and in general consider the mean squared error by choosing $q=2$.

\begin{figure}
    \centering
    \includegraphics[scale=0.40]{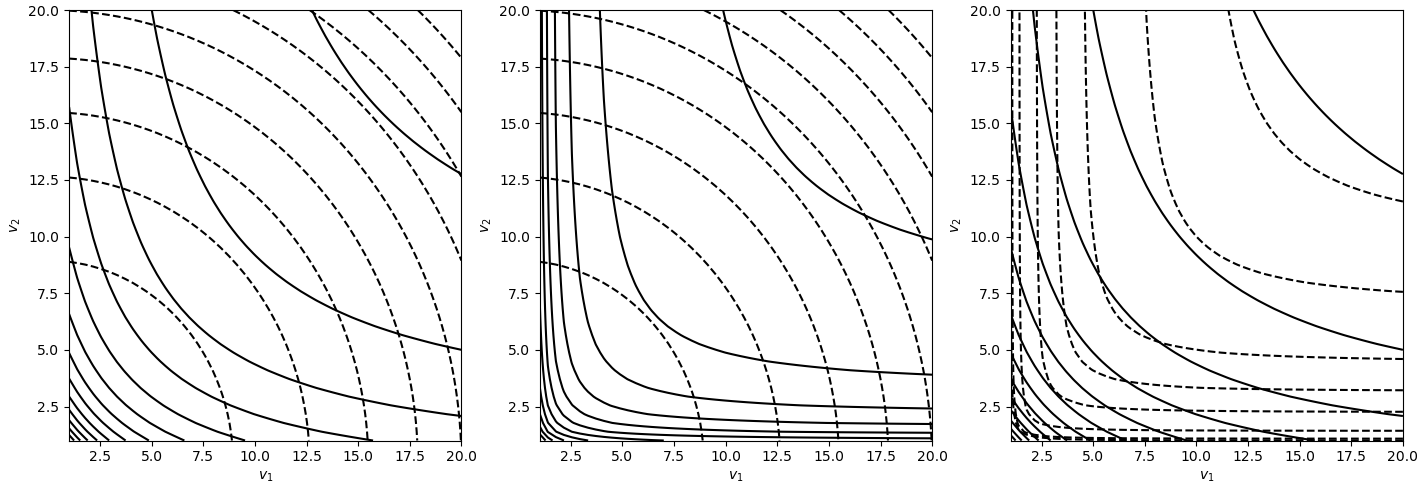}
    \caption{Sketch of the design problem~\eqref{eq:doe-incremental-reduced} for $n=2$ candidate points. Level lines of the objective $E(v)$ are drawn by solid lines, whereas those of the constraints are indicated by dashed lines. \emph{Left:} For $s>1$ there is a unique solution, which is likely not sparse. \emph{Middle:} A smaller correlation length $L<1$ makes sparsity even less likely. \emph{Right:} For $s<1$, the admissible sets are non-convex, and we may expect multiple local sparse minimizers.}
    \label{fig:design-sketch}
\end{figure}

In contrast, $s$ depends on the forward model $y$ and the simulation methodology. Looking at finite element simulations of order $r$ in $d$ dimensions, and assuming an optimal solver, we obtain $s = d/(2r)$. Consequently, for slowly converging linear finite elements in two or three space dimensions, a unique solution of the design problem exists, and we can expect it to be not sparse, see Fig.~\ref{fig:design-sketch} left. In contrast, quickly converging high order finite elements with $r\ge 2$ lead to $s<1$ and therefore non convex admissible sets. Their pronounced corners on the coordinate axes make the sparsity of a global minimizer likely, see Fig.~\ref{fig:design-sketch} right. This agrees with intuition: if increasing the accuracy at a specific sample point is computationally expensive, it is advantageous to distribute the work on a lower accuracy level to several points. On the other hand, if increasing the accuracy is relatively cheap, such as with quickly converging high-order finite elements, then it is often better to increase the accuracy of a single point, that, to some extent, shares its increased accuracy in a certain neighborhood.
While in the convex case $q\ge 2$ and $s\ge 1$ the optimization is straightforward with any nonlinear programming solver, the non convex case is more difficult. Fortunately, since the design problem is based on several relaxations and estimates, guaranteed global optimality is in practice not necessary. The expected sparsity structure suggests a particular heuristic approach: For $i=1,\dots,n+j$ consider $v=\underbar{v}+ae_i$ with $a>0$ such that $W(v)=\Delta W + W(\mathcal{D})$, i.e. the accuracy of only a single point is improved, and select the design $v$ with smallest objective. If this satisfies the necessary first-order conditions, accept it as solution. Otherwise, perform a local minimization starting from this point.

The algorithm proposed can also effectively operate in higher-dimensional parameter spaces. In order to add more additional points, we again approximate the maxima of the estimated error by equation \eqref{eq:utility} which is achieved by heuristically evaluating the error estimate on a grid of candidate points. However, in higher dimensions, this can be done by using (quasi) Monte Carlo sampling or multistart methods with local optimization.

\paragraph{Controlling the incremental budget.}
Different types of control can be used for the incremental budget $\Delta W$. It has been found that an exponential growth of the incremental budget is a suitable control, i.e we increase the incremental budget by 10$\%$ every iteration. In addition, we increase the budget by an additional 10$\%$ whenever the relative change in the estimated global error between iterations is too small.

\paragraph{Choice of kernel function.}
By assuming that the response surface is differentiable at any point and that the estimate of the derivative $f'$ is an integral component to determine the weighting factors, the choice of kernel functions is limited. Thus, all kernel functions of the form $\mathcal{K}(|x-x'|)$ are omitted.

\section{Numerical Examples}\label{sec:experiments}
In this chapter, we illustrate the presented concept numerically at two examples. In the first one, $y$ is given as an explicit arithmetic expression, while in the second example, it is given in form of a FEM simulation. 

\subsection{Analytical example}
We consider the rotated parabolic cylinder as model  $y$, i.e
\[
    y(p) = \left ( \cos(\phi)(p_1+p_2)+ \sin(\phi)(p_2-p_1) \right )^2 \quad \text{for $p \in \mathcal{X} = [0,1]^2$}, \; \phi \in \mathbb{R}_{> 0}
\]
We acquire $m=3$ measurements for $\phi \in \{0,2,4\}$. We further assume different accuracies of these independent measurements, so that the likelihood is of the form $\Sigma_{L} = 10^{-2}\mathrm{diag}(1,0.1,1)$ assuming anisotropic influence on $\tilde{w}$.

\begin{figure}[ht!]
    \begin{subfigure}[]{0.46\linewidth}
      \includegraphics[width=\textwidth]{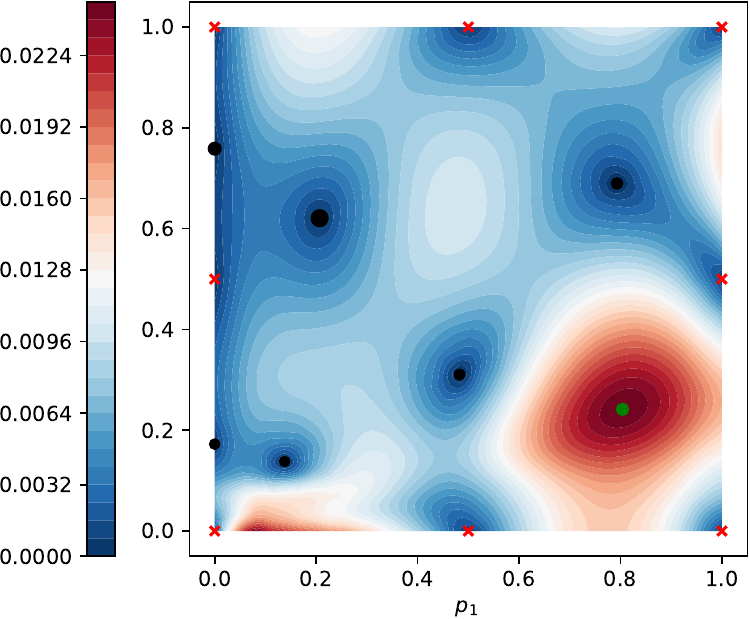}
      \label{subfig-1:smith}
    \end{subfigure}
    \hfill
    \begin{subfigure}[]{0.5\linewidth}
      \includegraphics[width=\textwidth]{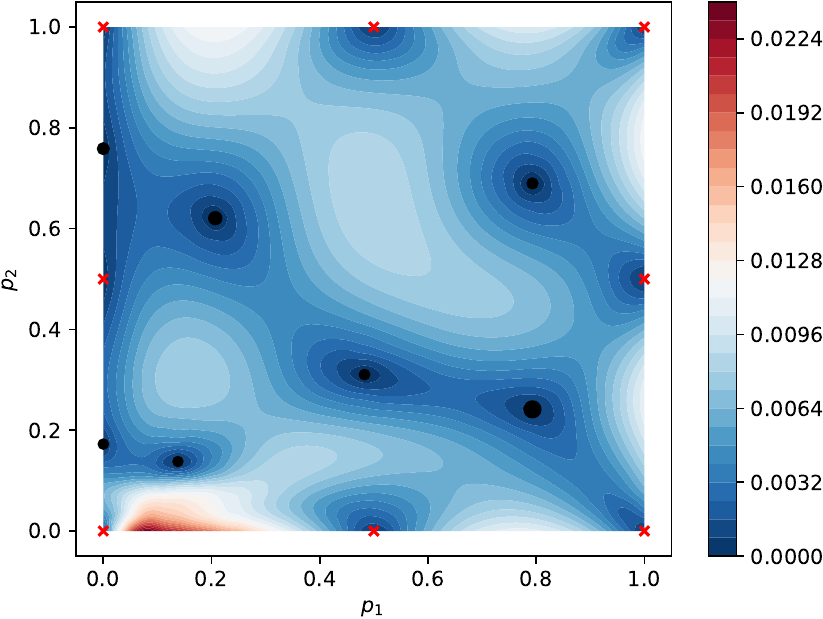}
      \label{subfig-2:smith}
    \end{subfigure}
\caption{Adaptively added data points are indicated via black dots, with size indicating accuracy -- small points indicate low accuracy and vice versa. Red crosses are initial data points. The color mapping shows the isolines of the estimated local reconstruction error evaluated on a dense grid of $10^3$ points. This design was obtained using an incremental budget of $\Delta W = 10^{4}$.
  \textit{Left:} Next to last design. The new point  $p=(0.804, 0.241)$ to be added is marked with a green point and indicates the maximum value of the acquisition function. \textit{Right:} Final parameter space after adaptive phase with $E_{\mathrm{ TOL }}(\mathcal{D})<\mathrm{TOL}$.}
\label{fig:2D_ParameterspaceAndError}
\end{figure}

For this example we start with an initial design of 8 evaluation points placed at the edges on $\partial\mathcal{X}$, see Fig.~\ref{fig:2D_ParameterspaceAndError}, with an evaluation variance of $\sigma^2 = 0.01$. We use the work model for quadratic finite elements in $d=2$ space dimensions, i.e. we set $s=1/2$.
Since we consider only a one low dimensional parameter space in this example, we omit Monte Carlo integration and use a standard quadrature to approximate the integral~\eqref{eq:MC-integration}. We use $N = 625$ points for the integral evaluation, arranged in an equidistant Cartesian grid. For the present example, this provides a sufficiently accurate approximation of the integral. 

Candidate points for inclusion into the domain are determined by evaluating the acquisition function~\eqref{eq:utility} on the grid mentioned above. In this example, only the point of maximum local reconstruction error is included in the design.
After successful minimization, the global error estimate is re-evaluated again.
The adaptive phase is terminated if the error falls below the desired tolerance, i.e $E_{\mathrm{ TOL }}(\mathcal{D}) \leq \mathrm{TOL}$. For this example, a global tolerance of $\mathrm{ TOL } \leq 10^{-2}$ is used.

We increase the incremental budget $\Delta W$ by $10\%$ at each iteration.

\paragraph{Results of adaptive phase.}

The result of the adaptive phase is shown in  Fig.~\ref{fig:2D_ParameterspaceAndError}, right. The left figure shows the previous design iteration. Black points represent adaptively added points. Note that the size of the symbol represents the accuracy of the data point -- small symbols mean low accuracy and vice versa. Color coding indicates the estimated local error. The local error reduction due to adding the design point at the position marked in green (Fig.~\ref{fig:2D_ParameterspaceAndError}, left)  is clearly visible. Note that the error at the point is not zero due to small but nonvanishing variance.

\begin{figure}[!ht]
    \begin{subfigure}[]{0.5\linewidth}
      \includegraphics[scale=0.52]{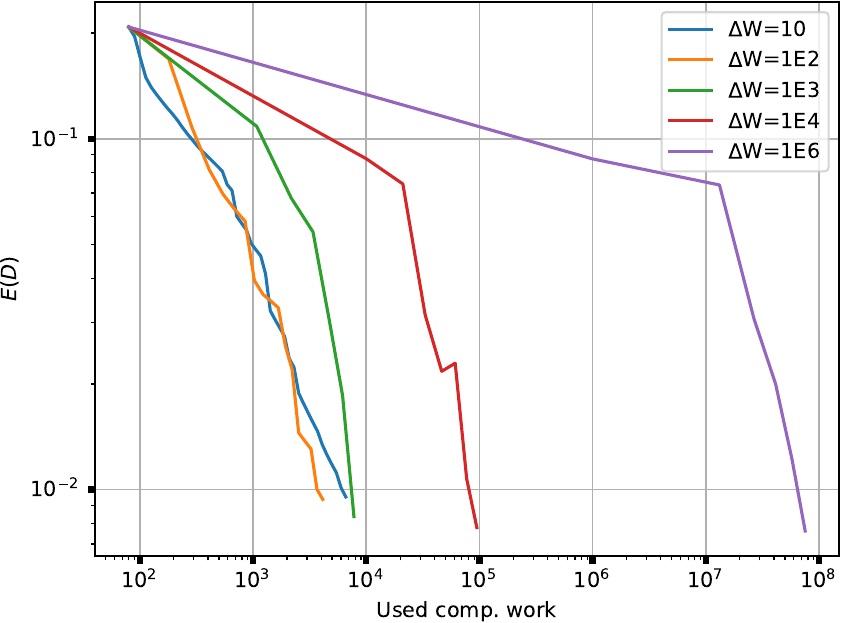}
    \end{subfigure}
    \hfill
    \begin{subfigure}[]{0.6\linewidth}
      \includegraphics[scale=0.52]{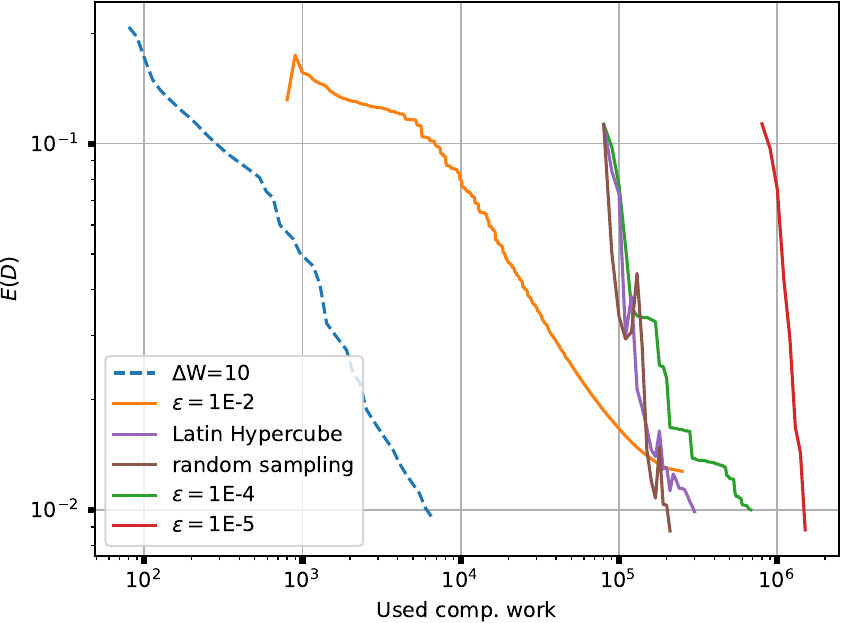}
    \end{subfigure}
    \caption{Estimated global error $E(\mathcal{D})$ versus computational work. \textit{Left}:  Different amounts of incremental work $\Delta W$. \textit{Right}: Different fixed evaluation accuracies $\epsilon$ compared with the curve for incremental work $\Delta W = 10$. Additionally we added two curves for a Latin Hypercube and a random sampling strategy. }
    \label{fig:2D_Errorplot}
\end{figure}

In Fig.~\ref{fig:2D_Errorplot}, left, the estimated global errors $E(\mathcal{D})$ are plotted against the total computational work $W$ for different incremental budgets $\Delta W$. 
The curve for the fully adaptive case (dashed line) and the semi-adaptive case (solid lines) is displayed in the right plot of Fig.~\ref{fig:2D_Errorplot}. In the semi-adaptive case, we utilize the acquisition function defined by equation  \eqref{eq:utility} for position optimization solely and subsequently adopt a constant evaluation error $\epsilon$ for the point $p$ determined in this manner.

By reducing $\Delta W$ on the left side, we can see that the curves  converge to a  minimal computational work necessary to achieve a given accuracy. Since here only one additional point is included into the design in each iteration, a rather small incremental budget yields the best performance.

The plots representing various incremental budgets $\Delta W$ are showcasing an initial decline in global error that intensifies as the incremental budget increases. This pattern arises from the fact, that the first adaptively added point attains a notably high evaluation accuracy when the budget is concentrated among only a few points. Notably, this effect becomes particularly pronounced in the curve associated with $\Delta W = 10^6$.

The magnitude of the subsequent trajectory slopes is heavily influenced by the size of the budget. For smaller budgets, the slopes are comparatively lower than those observed for larger budgets. This behavior aligns with expectations since a limited computational effort must be distributed across a larger number of points, resulting in a minimal increase in evaluation accuracy and consequently leading to a modest reduction in the estimate of global error. Conversely, for larger budgets, a substantial amount of computational work can be allocated among the points, resulting in a greater reduction of global error within a single iteration.

In the curve associated with $\Delta W = 10^4$, an abrupt surge in error occurs during one of the iterations. This occurrence can potentially be attributed to the utilization of suboptimal hyperparameters, resulting in a distorted derivative estimate and subsequently inaccurate calculation of the weighting factors $\tilde{w}$. As a consequence, the global error may increase as a result. To address this issue, employing multistart methods during hyperparameter optimization can effectively mitigate the impact of suboptimal settings and ensure the attainment of optimal parameter configurations.

To quantify the performance improvement of combined adaptivity in position and evaluation accuracy, the figure on the right shows the global error over computational work when using the fully adaptive algorithm for $\Delta W = 10$ (blue dashed line) and three different fixed evaluation errors $\epsilon$ (solid lines). In addition, we compare the fully adaptive case with a Latin hypercube (solid purple) and random sampling strategy (solid brown).

It is evident that the algorithm diverges for $\epsilon = 10^{-2}$ (orange) and fails to reach the desired accuracy. This is due to the fixed evaluation tolerance, which limits the surrogate model's accuracy to a maximum of $\epsilon$, regardless of the number of evaluation points used. It should be noted that although the horizontal progression could continue, we terminated the calculations prematurely because the error remained nearly constant.

If we lower the uncertainty to $\epsilon = 10^{-4}$ (green), the optimal case is (just) reached, i.e. the algorithm converges, the error curve takes on an almost horizontal course and thus the budget is optimally utilised. In the LH and random sampling settings, we also assigned an evaluation error of $\epsilon = 10^{-4}$ to the initial points. Both sampling strategies exhibit slightly faster convergence and demand less computational effort compared to the semi-adaptive approach. Nevertheless, when compared directly with the fully adaptive method, the amount of computational effort required is approximately 45 times higher.

On the far right the curve for $\epsilon =  10^{-5}$ is shown. It shows a rapid decrease of the error and thus rapid convergence, using approximately $W \approx  2 \cdot 10^6$. 
This behavior is expected since the semi-adaptive algorithm adds points with high evaluation accuracy in each iteration, which leads to a significant decrease in variance within the surrogate model and consequently a rapid reduction in the global error.

The fully adaptive algorithm outperforms the position adaptive algorithm in terms of computational efficiency by a factor 100. Full adaptivity requires a computational cost of $W \approx 8\cdot10^3$, while position adaptivity requires $W\approx 8\cdot 10^5$ for the optimal case. 

\paragraph{Reliability of local error estimator.}
The error model~\eqref{eq:globalerror_integral} is based on several assumptions, bounds and estimates, and may therefore not capture the actual error in identified parameters correctly. We therefore compare the estimated global error as formulated in the error model $E(\mathcal{D})$ to the actually obtained errors.

For 1600 points $p_i$, sampled randomly from $\mathcal{X}$, we compute the error estimate $\tilde e_i := w(p_i)\epsilon(p_i)$ on one hand, and obtain samples $e_{ik} := \|p(y(p_i)+\delta_{i,k})-p_\mathcal{D}(y(p_i)+\delta_{i,k})\|_2$, $k=1,\dots,K$ of the actual parameter deviation on the other hand, where $\delta_{i,k}$ are realizations of the measurement error distributed as $\mathcal{N}(0,\Sigma_l)$. Note that this involves minimizing the negative log likelihoods $J$ and $J_\mathcal{D}$ for the exact and the surrogate model, respectively.

We create two histograms, first for the estimated local errors $\tilde e_i$ and second for their sample mean $e_i := K^{-1} \sum_{k=1}^K e_{ik}$. 
\begin{figure}[!ht]
    \begin{subfigure}[]{0.5\linewidth}
        \includegraphics[scale=0.50]{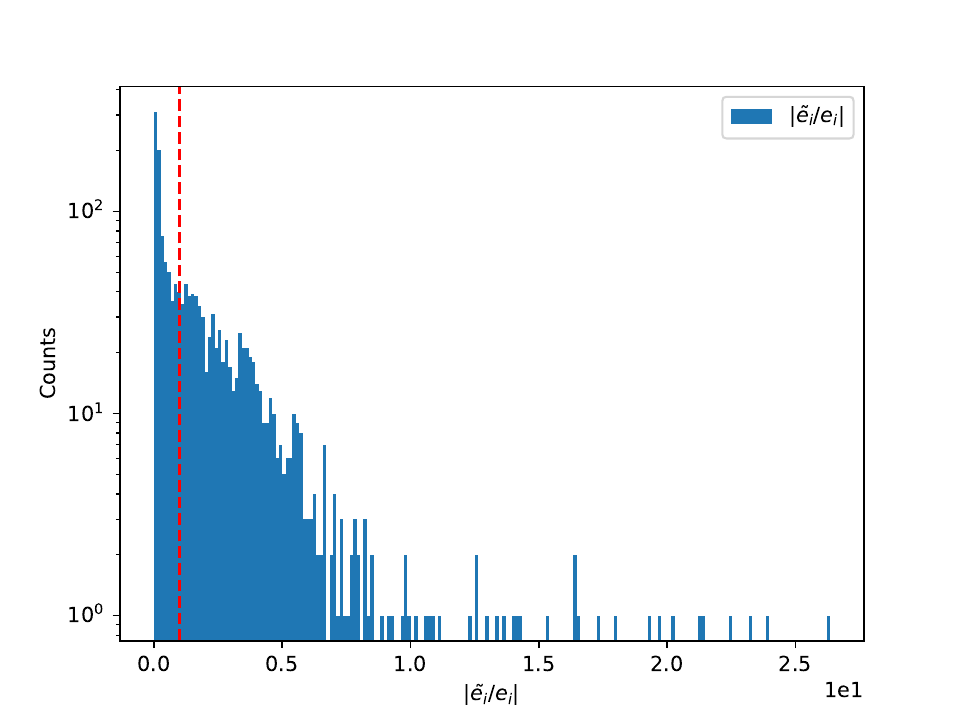}
    \end{subfigure}
        \begin{subfigure}[]{0.5\linewidth}
        \includegraphics[scale=0.50]{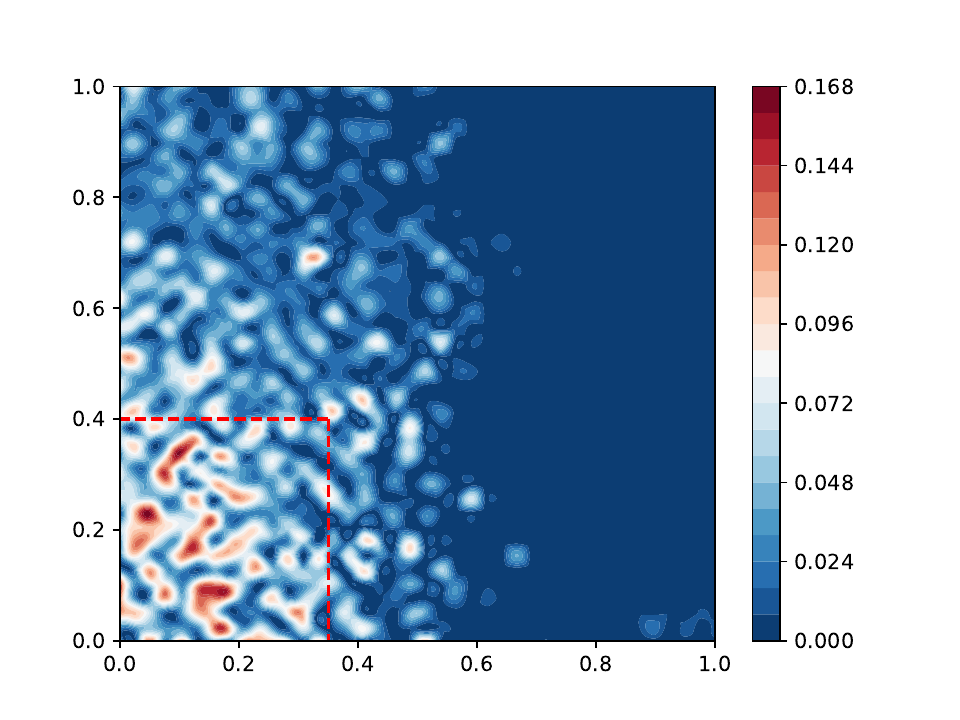}
    \end{subfigure}
    \caption{\textit{Left}: Log-histogram of $\tilde{e}\cdot{e}^{-1}$. \textit{Right}: Contour plot of $e_i$.}
    \label{fig:2D_Errorhist}
\end{figure}

The histogram Fig.~\ref{fig:2D_Errorhist} shows the error distribution for $\tilde{e}_i\cdot{e}_i^{-1}$ , i.e., the fraction of the estimated and the mean real parameter error. The red line represents the line $\tilde{e}_i\cdot{e}_i^{-1} = 1$, so all values to the left of it represent an underestimate and all values to the right of it represent an overestimate of the true reconstruction error. It shows that our assumptions are reasonable while underestimating the real reconstruction error. The large number of values on the right side, i.e. underestimated values, can be explained with the help of the right figure. The highest reconstruction errors are in the region enclosed by the red dashed lines (see. Fig.~\ref{fig:2D_Errorhist}, right), which was not sampled densely enough within the adaptive phase. By  construction, $y(p)$ changes its  value rapidly, which cannot be reproduced sufficiently well by the regression, which results in larger reconstruction errors.

The accuracy of the error estimate is dependent on two factors. Firstly, it relies on the accuracy of the GPR error estimate. Secondly, it is influenced by the linearization error of the error propagation estimate. Specifically, the accuracy of the GPR error estimate depends on the covariance model obtained through hyperparameter optimization.

Although it may be possible to demonstrate the reliability of the error estimator based on a priori assumptions regarding the regularity of the forward model and the spacing of evaluation points, this task is highly challenging and beyond the scope of this paper. Furthermore, even if a theoretical guarantee existed, verifying its assumptions in practical problems would be difficult.

Therefore, in this paper, we choose to focus on an empirical investigation of the reliability and efficiency of the error estimator, as we believe it is more practical and feasible at this stage.

\begin{figure}[!htb] 
  \centering
     \includegraphics[scale=0.6]{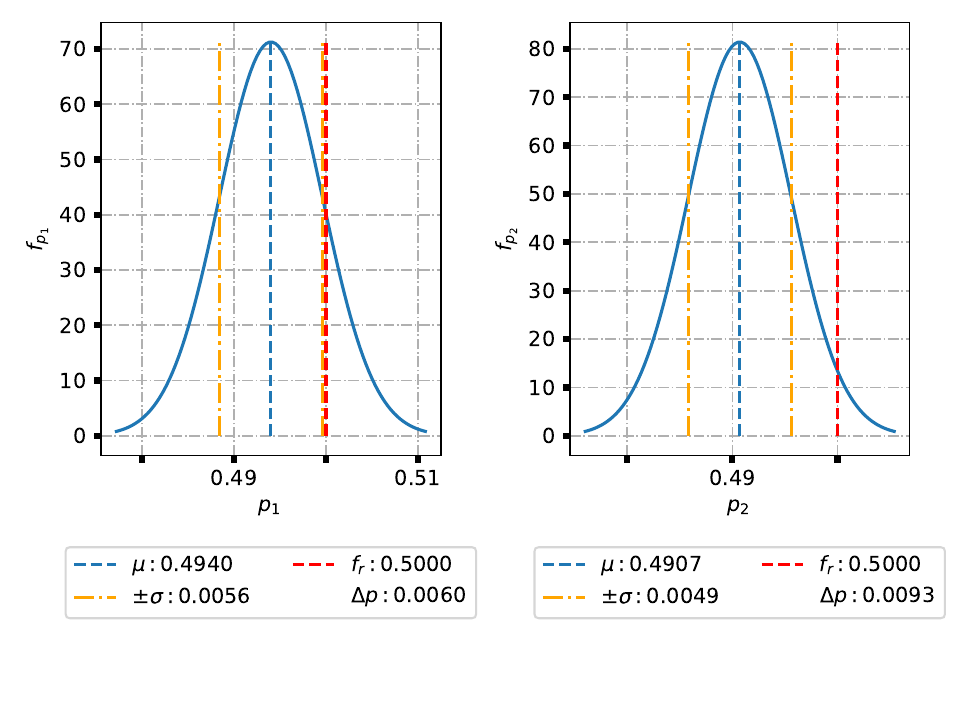}
  \caption{ Plot of the marginal distribution $f_{p_1}$ and $f_{p_2}$. Parameter reconstruction results for parameters $p_1$ and $p_2$ are shown by the red-blue dashed line. The artificial real parameters were chosen to be $p_{\textrm{real},1}=0.5$ and $p_{\textrm{real},2}=0.5$ and are shown by the red dashed line. The yellow dashed line indicates one standard deviation of the solution.}
 \label{fig:2D_reconstruction}
\end{figure}

\paragraph{Parameter reconstruction.}
Fig.~\ref{fig:2D_reconstruction} shows the marginal posterior distributions at the maximum likelihood point estimate for the reconstructed parameters given the adapted surrogate model (Fig.~\ref{fig:2D_ParameterspaceAndError}). With $p_{\textrm{recon.},1}=0.4940 \pm 0.0056$ and $p_{\textrm{recon.},2}=0.4907\pm 0.0049$, the results are consistent with the true parameters. Moreover, the absolute errors $\Delta p_{\textrm{recon.},1} = 6.0 \cdot 10^{-3}$ and $\Delta p_{\textrm{recon.},2} = 9.3 \cdot 10^{-4}$ are within the desired tolerance of $E_{\mathrm{ TOL }}(\mathcal{D})$.

Since we obtain the reconstructed parameter via a maximum likelihood point estimate, and the posterior is unimodal and rather narrow, we can employ the Laplace method of uncertainty quantification by using the diagonal elements of the covariance approximation 
\begin{equation}\label{eq:laplace}
    \sigma^2(p_{\text{MAP}}) := J''_{ii}(p_{\text{MAP}})^{-1} = \left((y'(p_{\text{MAP}}))^{T} (\Sigma_l+\Sigma_{\mathrm{GPR}})^{-1} y'(p_{\text{MAP}}) + \Sigma_p^{-1} \right)_{ii}^{-1}.
\end{equation}

Note that next to the likelihood covariance $\Sigma_l$ we also need to incorporate the surrogate variance $\Sigma_{\mathrm{GPR}}$ in order to include to surrogate model error.

\subsection{FEM Example}
In this example, we examine the behavior of the fully adaptive design construction using finite element simulations as forward model $y$. The algorithm has been worked out assuming the evaluation errors being identically, independently, and normally distributed. This does not hold for finite element discretization errors, which are correlated for close evaluation points $p$, and are not normally distributed due to the discrete nature of mesh refinement. However, the actual discretization error distribution is virtually impossible to describe reliably a priori, such that a coarse but deliberately simple Gaussian model as used here may nevertheless be a good choice. We therefore investigate, how the adaptive design construction copes in a setting not covered by the theoretical assumptions.

\begin{figure}[!htb] 
  \centering
     \includegraphics[scale=1.35]{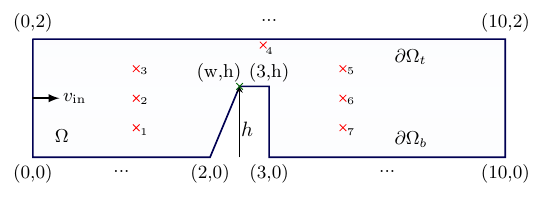}
  \caption{
  FEM example setup. The fluid flows into the domain $\Omega$  at the left boundary, while at the top and bottom boundary no-slip, i.e. $v=0$ boundary conditions are used. The flow velocities are measured at the red crosses. Both $w$ and $h$ are varied to generate the training data. Not that this sketch is not to scale.}
 \label{fig:2D_FEMsetup}
\end{figure}

We solve the Stokes equation 
\begin{align*}
    -\nu \Delta v+\nabla p&=0 && \text{in $\Omega$} \\ 
    \mathop\mathrm{div} v &=0 && \text{in $\Omega$} \\
    v&= v_0 &&\text{on $\partial\Omega_t \cup  \partial\Omega_b$ }, 
\end{align*}
with $\Omega = (\mathopen] 0, 10 \mathclose[ \times \mathopen] 0, 2 \mathclose[) \backslash \Omega_{\text{obs}}$ being a flow channel with an obstacle attached to a wall, see Fig.~\ref{fig:2D_FEMsetup}. The region of the obstacle $\Omega_{\text{obs}}$ is bounded by the convex polygon with corner points $(2,0), (w,h), (3,h)$ and $(3,0)$.
We consider homogeneous Dirichlet boundary conditions for the velocity at the top and bottom boundary, i.e. $v|_{\partial \Omega_t \cup \partial \Omega_b} = 0 $. To avoid edge singularities, we use a Hagen-Poiseuille profile for the incoming fluid on the left side and set the viscosity to $\nu = 1$. 

We consider the problem to be parameterized by two parameters defining the obstacle. The first parameter 
$w\in[2.1,2.9]$ describes the horizontal position of the obstacle, while the second parameter $h\in[1.5,1.9]$  specifies its height.
The velocity is measured at points marked by the red crosses, such that we consider $m=7$. Note that the middle evaluation point is located close to the top boundary  and inside the computational domain $\Omega$ also for large values of $h$. The problem is discretized with Taylor-Hood elements in the finite element toolbox Kaskade 7~\cite{GoetschelWeiserSchiela2012}. Since we are dealing with a two dimensional problem and therefore quadratic finite elements regarding the velocity discretization, we set $s=1/2$ within the work model.

As initial standard deviation of all data points we use the discretization error provided by Kaskade's internal hierarchical error estimator based on the difference of the computed FE solution and its projection onto a lower order ansatz space. In order to create a finer mesh, we use uniform mesh refinement over the whole domain $\Omega$.
We again choose $\mathrm{TOL} \leq 10^{-2}$ as the upper bound on the global error and also choose the likelihood being diagonal with $\Sigma_L = \mathrm{diag}(10^{-3}_1,\dots, 10_7^{-3})$, i.e the same accuracy for every measurement point and therefore isotropic influence on $\Tilde{w}$.

In contrast to the analytical example, the actual evaluation accuracy of the simulations may be higher than specified in the optimization, due to the discrete levels of refinement at which the simulations can be performed. The difference in accuracy within a design iteration must be tracked and the (incremental) computational budget for the following design iteration adjusted accordingly. 
If the evaluation incurs a higher computational effort than the incremental budget $\Delta W$ allocated for the optimization iteration, the difference $W(v(p_i))-\Delta W > 0$ is no longer available for the next design iteration.

If the optimization leads only to a small change in the evaluation accuracy of a data point, we refrain from translating this change, since only significant changes can be represented by a next higher discretization level.

\begin{figure}[htb!] 
  \centering
     \includegraphics[scale=0.6]{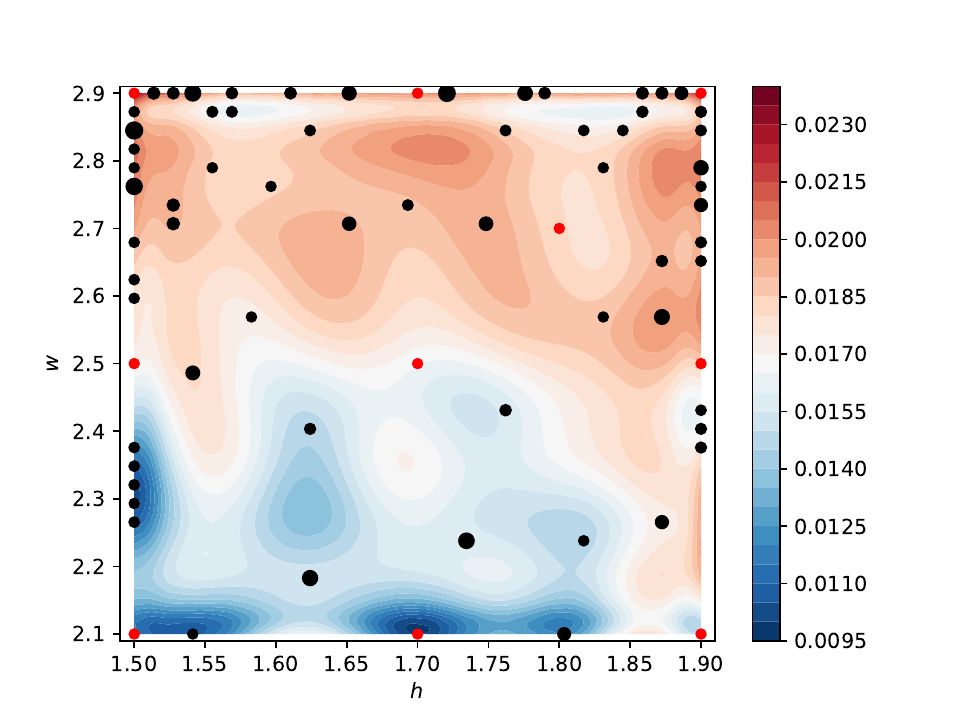}
  \caption{Parameter space after the adaptive phase. Adaptively added data points are indicated via black dots, with size indicating accuracy -- small points indicate low accuracy and vice versa. Red points are initial data points. The color mapping shows the isolines of the estimated local reconstruction error. }
 \label{fig:2DFEM_parameterspace}
\end{figure}

\paragraph{Adaptive phase.}
Fig.~\ref{fig:2DFEM_parameterspace} shows the local reconstruction error after the adaptive phase for $\Delta W_0 = 10^{2}$. As expected, more accurate data points, again indicated by point size, were added to the parameter space in the region with large $h$ and thus tight constrictions. This is due to the fact that there the velocity changes particularly strongly with small changes in the height of the constriction and thus strongly affects the weighting factors $\Tilde{w}$. 

\begin{figure}[htb!] 
  \centering
     \includegraphics[scale=0.60]{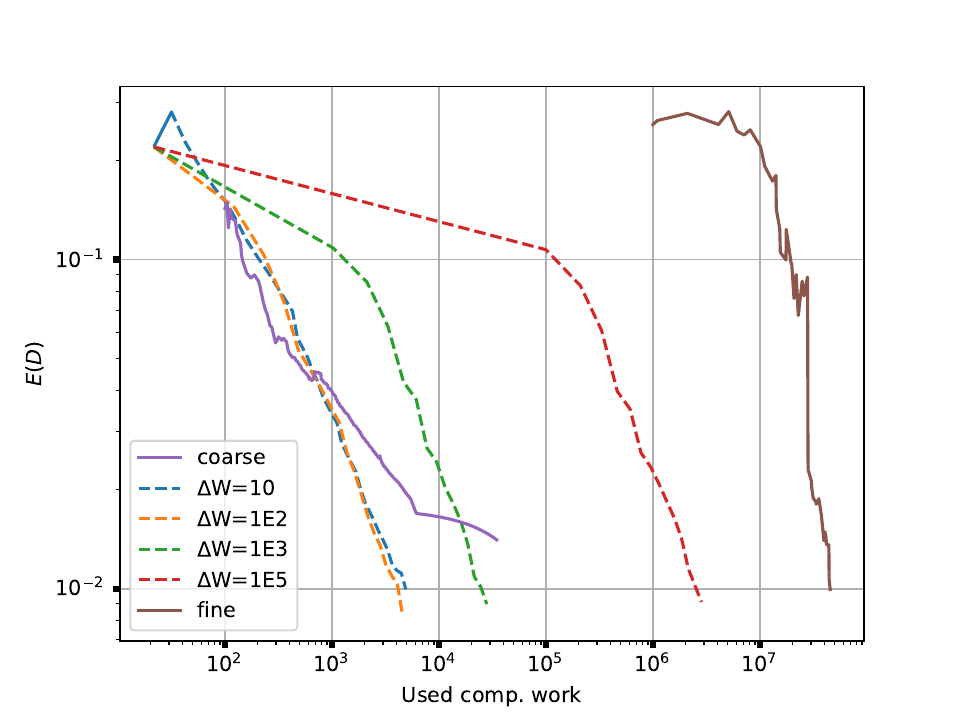}
      \caption{Estimated global error $E(\mathcal{D})$ versus computational work. The fully adaptive curves are shown as dashed lines, the curves for the position-adapted algorithm are shown as solid lines, and the curves for the sequential designs are shown as dashed lines.}
 \label{fig:2DFEM_performanceboost}
\end{figure}

\paragraph{Performance boost.}
Analogous to the analytic example, we examine the performance gain compared to the position-adaptive approach. Simulations are performed at the finest and coarsest (uniform) refinement levels, with the coarsest refinement level chosen such that the physical values are still reasonable. Fig.~\ref{fig:2DFEM_performanceboost} shows the estimated global error $E(\mathcal{D})$ plotted against the total computational effort.

The dashed lines represent the fully adaptive curves, while the solid lines depict the position-adapted curves. Similar to the analytical example, we observe a lower bound on the required computational effort, approximately $W\approx 6\cdot10^3$, as the incremental budget $\Delta W$ approaches zero. Among these curves, the purple curve represents the progression of the position-adaptive algorithm. Initially, it follows a linear trend, but after a certain number of iterations, it transforms into an almost horizontal curve. Consequently, the accuracy of the surrogate model experiences minimal improvement, and the desired accuracy remains unattained. However, by increasing the level of refinement, the brown curve is obtained. This curve converges to the desired accuracy, but it necessitates a significantly higher computational effort of approximately $W \approx 6\cdot 10^7$, which is approximately a thousand times greater than that of the fully adaptive approach.

\newpage
\paragraph{Reliability of local error estimator.}
Following the evaluation method of the analytical example, we again create a histogram for the fraction $\tilde{e}_i\cdot{e}_i^{-1}$, using $N=50$ this time and a contour plot of $e_{i}$ over the whole parameter space. 
\begin{figure}[!ht]
    \begin{subfigure}[]{0.5\linewidth}
        \includegraphics[scale=0.50]{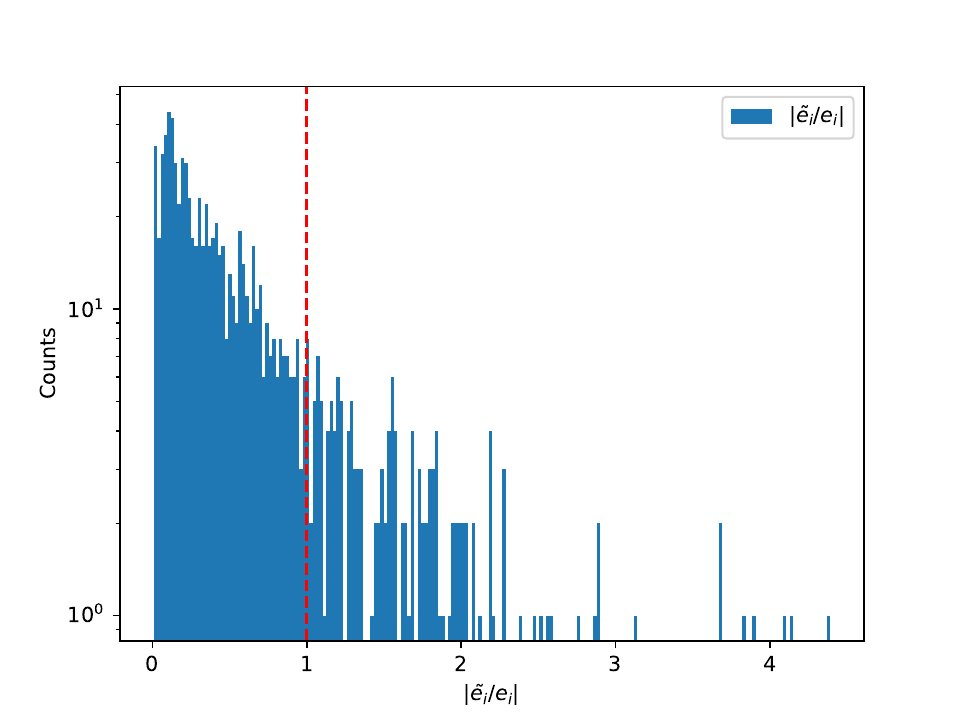}
    \end{subfigure}
        \begin{subfigure}[]{0.5\linewidth}
        \includegraphics[scale=0.50]{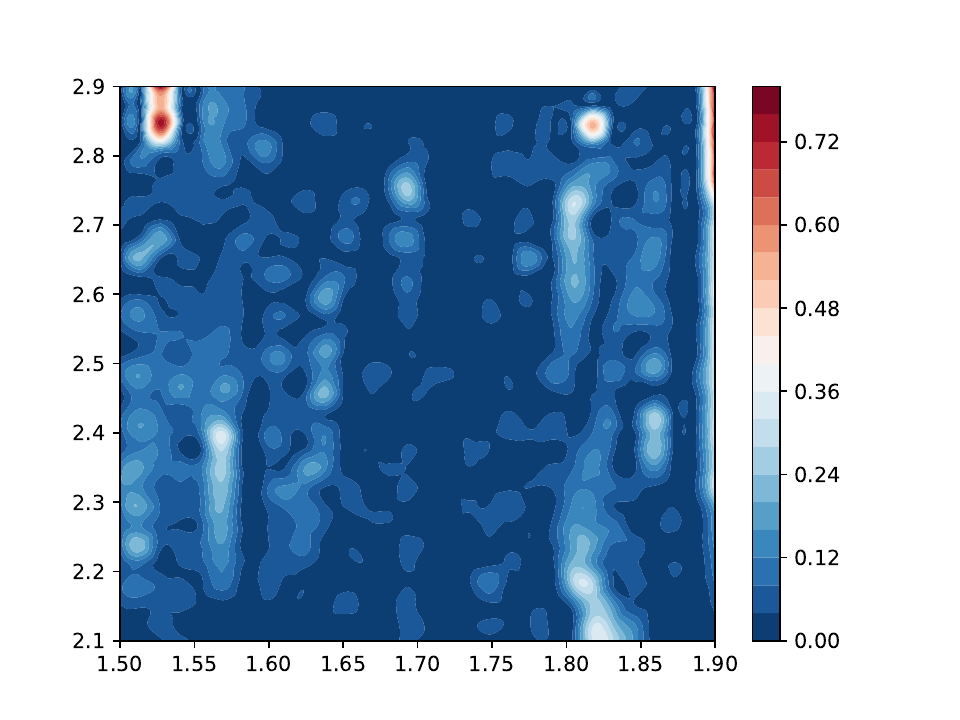}
    \end{subfigure}
    \caption{\textit{Left}: Log-histogram of $\tilde{e}_i\cdot{e}_i^{-1}$. \textit{Right}: Contour plot of $e_i$.}
    \label{fig:2D_FEM_Reliability}
\end{figure}
For the evaluation, we only consider the range up to $\tilde{e}_i \cdot e_i^{-1}\approx 2$, as we characterise the errors above this limit as outliers.  The histogram shows that we underestimate the error for the most part and can only reproduce the mean reconstruction error well in a small range. If we consult the figure on the right for an explanation, we see that these errors occur where the parameter space was not sampled densely enough or the data points have only a low evaluation accuracy. This leads to the conclusion that the assumptions made are too simple for real FEM data.

\paragraph{Parameter reconstruction.}
For the parameter reconstruction, we assume the parameter pair $p_{\mathrm{true}} = (1.6, 2.7)$ as true parameters and use FEM simulations at the highest discretization level in order to create the vector of experimental data $y_m$. 
\begin{figure}[!htb] 
  \centering
     \includegraphics[scale=0.6]{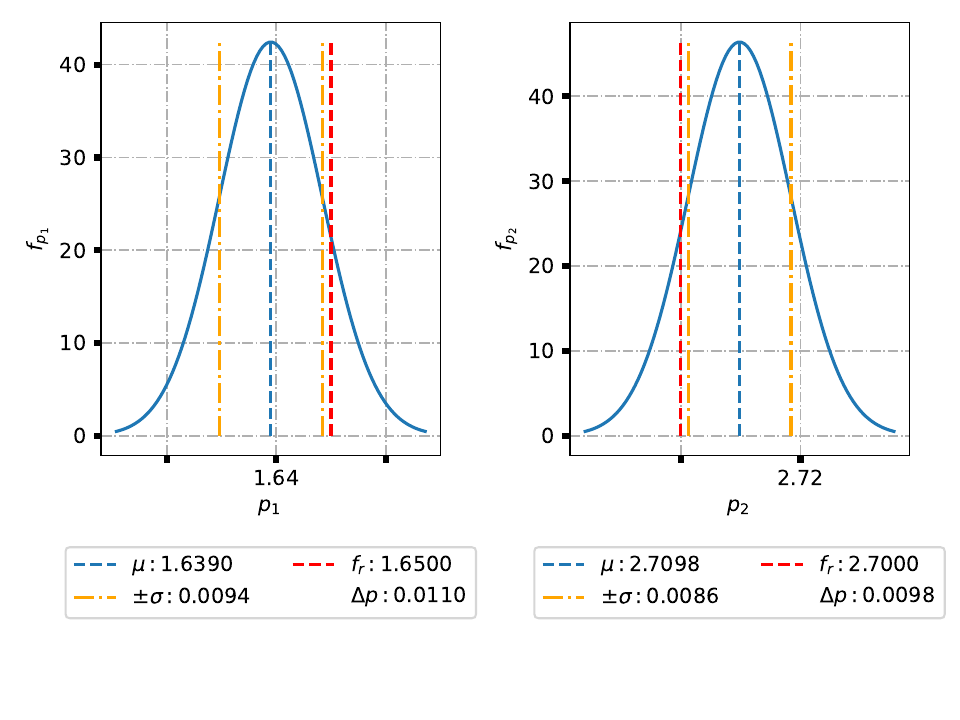}
  \caption{ Plot of the marginal distribution $f_{p_1}$ and $f_{p_2}$.
 Parameter reconstruction results for parameters $p_1$ and $p_2$ are shown by the red-blue dashed line. The artificial real parameters were chosen to be $p_{\textrm{true},1}=1.6$ and $p_{\textrm{true},2}=2.7$ and are shown by the red dashed line. The yellow dashed line indicates one standard deviation of the solution.}
 \label{fig:2DFEM_solution}
\end{figure}
Fig.~\ref{fig:2DFEM_solution} shows the maximum a posteriori estimates after successful reconstruction. The errors of the reconstructed parameters $\Delta p_1 = 3.6\cdot 10^{-3}$ and $\Delta p_2 = 4.3\cdot 10^{-3}$ are within desired tolerance of $E(\mathcal{D}) \leq 10^{-2}$ and are consistent with the true parameters. Again, in order to quantify the uncertainty of the reconstruction we use Laplace method ~\eqref{eq:laplace}.

\newpage
\section{Conclusion}
This study presents an adaptive algorithm for the optimal construction of surrogate models for use in parameter reconstruction of inverse problems. It provides significant savings in computational cost compared to a priori designs based on random parameter positions, low-discrepancy series, or position-adaptive strategies. An reduction of computational effort by a factor of $100-1000$ compared to a position adaptive scheme was observed in examples.

The reliability analysis shows that the assumptions made on the evaluation error distribution give good results for artificial data, but lead to an underestimation of the actual error for more complex FEM data. Nevertheless, the regression models were able to reconstruct the parameters with the desired accuracy for both examples.

The choice of error model and the kernel are kept simple in the first approach and can be built up more complicated and problem-adapted, i.e. a non Gaussian error model within the likelihood.

Numerical investigations have shown that the estimate of the derivative $f'$ depends very sensitively on the choice of the hyperparameter $L$, as well as strongly on the number of points within the system and their evaluation accuracy. Thus, the determination of the error weighting factors $w$ changes significantly at the initial adaptive phase and can thus lead to sub optimal positioning of new data points. More theoretical and  numerical investigations are necessary.

\paragraph{Funding.}
This work has been supported by Bundesministerium für Bildung und Forschung -- BMBF, project number 05M20ZAA (siMLopt).

% \printbibliography
\bibliographystyle{plain}
\bibliography{AdaptiveGPR}

\end{document}